%% file: optdir23.tex
\documentclass[final]{siamltex}
\usepackage{amsfonts,amsmath,mathrsfs}
\usepackage{amssymb}
\usepackage{hyperref}
\usepackage{graphicx}
\usepackage{wrapfig}
\usepackage{subfig}
\usepackage{float}
\usepackage{bbm}
\usepackage{algorithmic}
\usepackage{gensymb,color}
\usepackage{tikz}

\newcommand{\N}{\mathbb{N}}
\newcommand{\R}{\mathbb{R}}

\newcommand{\dx}{{\rm d}x}
\newcommand{\dy}{{\rm d}y}
\newtheorem{remark}[theorem]{Remark}
\newtheorem{algorithm}{Algorithm}

\DeclareMathOperator{\tr}{\rm tr}
\title{Sequentially optimized projections in X-ray imaging}

\author{M. Burger\footnotemark[2]
\and  A. Hauptmann\footnotemark[3]
\and T. Helin\footnotemark[4]
\and N. Hyv\"onen\footnotemark[5]
\and J.-P. Puska\footnotemark[5]
}

\begin{document}
\maketitle

\renewcommand{\thefootnote}{\fnsymbol{footnote}}
\footnotetext[2]{Friedrich-Alexander Universit\"at Erlangen-N\"urnberg, Department Mathematik, Cauerstrasse 11, 91058 Erlangen, Germany (martin.burger@fau.de). The work of MB  was supported by the Bundesministerium f\"ur Bildung und Forschung under the project id 05M16PMB (MED4D).}
\footnotetext[3]{University of Oulu, Research Unit of Mathematical Sciences, Oulu, Finland; University College London, Department of Computer Science, London, United Kingdom (andreas.hauptmann@oulu.fi). The work of AH was supported by the Academy of Finland (decision 312123).}
\footnotetext[4]{LUT University, School of Engineering Science, P.O.~Box 20, FI-53851 Lappeenranta, Finland. The work of TH was supported by the the Academy of Finland (decisions 320082 and 326961).}
\footnotetext[5]{Aalto University, Department of Mathematics and Systems Analysis, P.O.~Box 11100, FI-00076 Aalto, Finland (nuutti.hyvonen@aalto.fi, juha-pekka.puska@aalto.fi). The work of NH and JP was supported by the Academy of Finland (decision 312124).}

\begin{abstract}
  This work applies Bayesian experimental design to selecting optimal projection geometries in (discretized) parallel beam X-ray tomography assuming the prior and the additive noise are Gaussian. The introduced greedy exhaustive optimization algorithm proceeds sequentially, with the posterior distribution corresponding to the previous projections serving as the prior for determining the design parameters,~i.e.~the imaging angle and the lateral position of the source-receiver pair, for the next one. The algorithm allows redefining the region of interest after each projection as well as adapting parameters in the (original) prior to the measured data. Both A and D-optimality are considered, with emphasis on efficient evaluation of the corresponding objective functions. Two-dimensional numerical experiments demonstrate the functionality of the approach.
\end{abstract}

\renewcommand{\thefootnote}{\arabic{footnote}}



\pagestyle{myheadings}
\thispagestyle{plain}
\markboth{A. BURGER, A. HAUPTMANN, T. HELIN, N. HYV\"ONEN, AND J.-P. PUSKA}{OPTIMIZED PROJECTIONS IN X-RAY IMAGING}

\section{Introduction}
\label{sec:introduction}

Tomographic image reconstruction is a classical inverse problem with its history spanning over 100 years since the seminal work of Radon \cite{radon20051}. It gives rise to a wealth of applications in modern imaging for which research remains active to date. The contemporary mathematical questions are often driven by practical challenges such as limited data --- consider,~e.g.,~having a limited number of imaging angles or a limited field-of-view. There exists a number of proposed methods for image reconstruction given such incomplete data, however, much less is known about how to optimally design the imaging configuration if there are constraints on,~e.g.,~the number of angles, if individual projections do not cover all of the imaged object, or if there is a particular {\em region of interest} (ROI) that may vary as a function of time. Such design is desirable in any application where the exposure to radiation must be minimized or the acquisition of data is otherwise expensive.

Optimal experimental design (OED) is also an extensive field of research with a massive body of literature; for a general introduction we refer to \cite{atkinson2007optimum}. In particular, OED in inverse problems has taken leaps during the last decade due to the increase in computational resources. In the present paper, we adopt the Bayesian paradigm to inverse problems, that is, the observed data is used to update the prior knowledge about the unknown into a posterior probability distribution \cite{kaipio2006statistical, stuart_2010}. OED in Bayesian inference is formulated as minimization of a Bayesian cost with respect to the design parameters \cite{chaloner1995bayesian}; different cost functions correspond to different optimality criteria. Here we focus on the two most common approaches: \emph{A-optimality}, where a typical quadratic loss is considered, and \emph{D-optimality}, which maximizes the Kullback--Leibler divergence of the posterior distribution with respect to the prior.

For general high-dimensional inference problems, the minimization of the Bayes cost requires a huge computational effort since just evaluating the cost corresponds to (numerically) integrating the considered cost function against the joint probability distribution of the unknown and the data. Considerable effort has been put into formulating approximative schemes to solve the corresponding optimization task; in the framework of inverse problems this work is discussed below. However, for linear Bayesian inference problems with Gaussian prior and likelihood, it is well-known that the integral defining the Bayes cost for A or D-optimality can be evaluated analytically as the trace or the determinant of the posterior covariance, respectively, leading to more attractive optimization problems~\cite{chaloner1995bayesian}. This is the general setting considered in this work.

\subsection{Our contribution}
Our paper focuses on discussing benefits of OED in linear X-ray imaging. Our contributions to this subject are as follows:
\begin{itemize}
\item We propose a greedy exhaustive sequential method for optimizing the imaging configuration in X-ray tomography in the spirit of \cite{fohring2016adaptive} and discuss its efficient numerical implementation. More precisely, our work proposes how to optimize the next imaging angle and the lateral position of the source-receiver pair if the geometric specifications of the previous projections are known. At each step, the posterior covariance corresponding to the previous projections is fed as the prior covariance for the optimization of the next imaging angle and the lateral position. We discuss both A and D-optimal designs. As the dimension of the new measurement at each step is considerably lower than that of the unknown, the computational cost of evaluating the necessary posterior traces and determinants can be significantly lowered by resorting to the Woodbury matrix identity and the matrix determinant lemma without employing Monte Carlo estimators (cf.~\cite{alexanderian2016fast,alexanderian2018efficient,avron2011randomized}).
\item If the covariance matrices for the (original) Gaussian prior and the additive Gaussian noise process are known, the sequentially optimal projections can be determined before performing any actual measurements. However, we also introduce an optimization algorithm based on a varying ROI. Indeed, the ROI can be adapted to the reconstruction after a new projection image becomes available either algorithmically or based on expert knowledge, and subsequently the next optimal projection can be determined accounting for the change in the ROI. A related idea for the A-optimality criterion is suggested in~\cite{fohring2016adaptive}.
\item We develop a data-driven variant of our optimization method that simultaneously adapts the (original) prior distribution to the data observed at each step. As a numerical example, we demonstrate that the method is able to identify the true correlation length used in generating an imaged target.
\end{itemize}

It should be emphasized that the forms of the optimization target functions corresponding to both A and D-optimality criteria are well-known \cite{chaloner1995bayesian}. However, linear X-ray imaging provides an interesting example for further analysis of the setting, not only due to its high-dimensionality, but also since the design space (i.e.~the projection geometry) is continuous by nature, which leads to challenging computational considerations. We adopt a sequential approach to the optimization due to its computational tractability and in order to include feedback from the observed data in the design process. Notice that sequential optimization is in general suboptimal if the total number of projections is fixed {\em a priori} and no feedback from the measured data is included in the process~\cite{huan2016sequential}.

\subsection{Literature review}

In recent years, OED for large-scale inverse problems has gained substantial attention; see,~e.g.,~\cite{alexanderian2014optimal, alexanderian2016fast, alexanderian2018efficient, attia2018goal, haber2008numerical, haber2009numerical, haber2012numerical, Hannukainen20,hyvonen2014eit, khodja2010guided,long2015fast,long2013fast}. In particular, we acknowledge that the study of A and D-optimal designs has been extended to infinite-dimensional Hilbert spaces in \cite{alexanderian2016bayesian}. The computational feasibility of experimental design in large-scale problems is often based on linearization and the use of Monte Carlo trace estimators such as the ones developed in \cite{avron2011randomized}. However, there is also an interesting avenue of research presented in \cite{huan2014gradient,huan2013simulation, huan2016sequential}, where the authors develop design methods based on dynamical programming.

Optimization of projection angles in X-ray imaging has previously been considered in \cite{ruthotto2017optimal}. The authors explore two interesting approaches for empirical A-optimal design in constrained problems based on training data: On the one hand, they solve a sparse design out of a predetermined set of angles, which is an approach aligned with previous literature. On the other hand, they propose a gradient-based optimization scheme. However, there seems to be no guarantee of being able to find a global optimum and the presented numerical experiments are limited to two angles. The present paper is closely related to \cite{ruthotto2017optimal} as our approach is also based on a prior. As novelties compared to~\cite{ruthotto2017optimal}, we consider the lateral position of the source-receiver pair as a second design parameter, our sequential design enables incorporating feedback from the observed data in the process, and the exhaustive optimization algorithm facilitates circumventing difficulties arising from the nonconvexity of the objective functions. In fact, our numerical examples clearly demonstrate that the optimization of the projection geometry in X-ray imaging suffers from the existence of local optima.

In terms of adaptivity, our work draws inspiration from \cite{fohring2016adaptive}, where A-optimal sequential design for dynamical inverse problems is considered. Our setup can be derived from \cite{fohring2016adaptive} by assuming static noise-free dynamics. However, in \cite{fohring2016adaptive} the adaptation of design is obtained by introducing a so-called monitor function that measures the difference in the last two reconstructions. In our framework, the approach of \cite{fohring2016adaptive} would roughly correspond to choosing the ROI to be an area where such a difference is high.

From an application's point of view, our proposed method fits in the larger field of sparse tomography \cite{hamalainen2013sparse,Siltanen03}, but moreover to applications with strict dose limitations and obstructions to the imaging region \cite{borg2017reduction}. Such limitations might be due to safety reasons, experimental environment, or even expensive beamtime at synchrotron facilities \cite{Batenburg2018advanced}. For instance, in image-guided radiotherapy angles are chosen to radiate only a cancerous region of interest, without damaging essential surrounding tissue, using the knowledge obtained during the measurements,~cf.~\cite{jaffray2012image}. Here, the proposed OED framework could provide insights into treatment planning and creation of automated algorithms to exclude critical regions from the beam direction.

This paper is organized as follows. In Section \ref{sec:mm}, we introduce our discretized parallel beam framework for two-dimensional X-ray imaging. Section~\ref{sec:bayes} reviews the basics of finite-dimensional Bayesian OED and considers efficient numerical evaluation of the objective functions corresponding to A and D-optimality. The greedy exhaustive optimization algorithm and its adaptive variants are presented in Section~\ref{sec:opti}. Finally, the numerical experiments are documented in Section~\ref{sec:numerics}, and some conclusions are drawn in Section~\ref{sec:concl}. For completeness and to facilitate reading the paper, an appendix presents brief derivations of the target functions for A and D-optimal designs.

\section{Measurement model and its discretization}
\label{sec:mm}
The measurements are modeled as parallel beam tomography, where multiple parallel rays are directed into the object $D \subset \R^2$, and the resulting intensities of the rays are measured at detectors~\cite{Natterer01}. The attenuation is described by the equation
\begin{equation}
	\label{eq:model}
	I = I_0 \exp \left( - \int_{L} u \, {\rm d}s \right),
\end{equation}
where $L$ is the line along which the considered ray travels, $I_0$ is the intensity of the X-ray before entering the object and $u: D \to \R_+$ is the absorption. Obviously, \eqref{eq:model} can equivalently be given as
$$
\log(I_0) - \log(I) = \int_{L} u(x) \, {\rm d} s.
$$
In particular, the difference between the logarithms of the  emitted and measured intensities is typically considered as the available data when X-ray tomography is tackled mathematically.

After discretizing the imaged domain into pixels, the forward operator, mapping the discretized absorption to a {\em single} set of log-intensity measurements at the detectors can be approximated by a matrix $R \in \mathbb{R}^{m \times n}$, where $m$ and $n$ are the  number of detectors and the number of pixels, respectively (see,~e.g.,~\cite{Siltanen03}); typically the dimension of the unknown is higher than the number of pixels in a {\em single} projection image, i.e.~$m\ll n$.  Each ray is parametrized by its (signed) distance $s$ to the origin and its angle $\phi$, say, relative to the positive horizontal axis. For a single projection image with $m$ detectors, one thus needs to specify $m+1$ parameters $(s_1, \ldots, s_m ,\phi)$. We assume the rays/detectors to be equally spaced, but in our setting the complete parallel beam source-receiver pair does not typically cover the whole domain, that is, a single measurement with a given angle only produces a projection image of a part of $D$.

Take note that $R$ is very sparse, which makes corresponding matrix multiplications relatively inexpensive computationally. Moreover, at least if $m \lesssim \sqrt{n}$, it is reasonable to expect that $R$ has a full (effective) rank, i.e.,~${\rm rank}(R) = m$, with all its singular values being of the same order of magnitude. In particular, no projection matrix considered in our numerical experiments has a Euclidean norm based condition number that is higher than $2$.

We assume the domain $D$ consists of three distinct regions:
\begin{itemize}
	\item {\em ROI}: the area about which we want to recover information,
\item \textit{Obstruction}: a nuisance region that obstructs the propagation of X-rays,
\item \textit{Background}: the rest of the object.
\end{itemize}
By resorting to this kind of a division, one can model settings where only a certain part of the target object is of interest. Such situations could arise,~e.g.,~in a medical application if one has already determined the location of  an organ, tumor or other abnormality that is to be further monitored. Moreover, in an industrial application one may only be interested in a certain component or section of some structure. In addition, the introduced division enables modeling obstructions that are difficult to image through, such as thick bone structures or other highly attenuating materials.

\section{Bayesian inversion and experimental design}
\label{sec:bayes}
In Bayesian inversion all parameters carrying uncertainty are treated as random variables \cite{kaipio2006statistical,stuart_2010}. The prior probability distributions for these parameters reflect the available information before the measurements are carried out. A measurement is modeled as a realization of a random variable depending on both the noise process and the random parameters in the forward model, as well as on the so-called design parameters that are deterministic and define the measurement setup. The Bayes’ formula is then employed to form the posterior probability density that updates the prior based on the information in the measurement. In our setting, the design variables are the projection angle and the lateral position of the parallel beam source-receiver pair, and our ultimate aim is to (sequentially) choose these parameters according to either the A or D-optimality criterion of Bayesian experimental design \cite{alexanderian2016bayesian,chaloner1995bayesian}.

In this section, we only consider one step of the optimization process, that is, we assume to have a prior covariance (or its inverse) at hand and we consider optimally choosing the design parameters that define a single parallel beam projection. The computational efficiency is in the focus of our attention. The actual sequential optimization algorithm is presented in Section~\ref{sec:opti}.

\subsection{Bayes' formula}

Let $y \in \R^m$ represent the (noisy) projection image and $p \in \R^2$ be the design variable, and suppose our prior information on the (discretized) absorption distribution is encoded in a probability density $\pi_{\rm pr}: \R^{n} \to \R_+$. By Bayes' formula, the posterior density for the (randomized) absorption $X$ reads
\begin{equation}
  \label{eq:Bayes}
\pi(x \, | \, y;  p) = \frac{\pi(y \, | \, x;  p) \pi_{\rm pr}(x)}{\pi(y; p)},  \qquad x \in \R^{n},
\end{equation}
where $\pi(y \, | \, \cdot \, ;  p): \R^{n} \to \R_+$ is the so-called likelihood function and the normalizing term in the denominator is the marginal density of the measurement $Y$ evaluated at the available data $y$.

In this work we assume the prior is Gaussian, i.e.~$X \sim \mathcal{N}(x_{\rm pr}, \Gamma_{\rm pr})$, and the measurement can be modeled as a realization of the random variable
\begin{equation}
  \label{eq:meas_model}
Y = R(p) X + N,
\end{equation}
where $N \sim \mathcal{N}(0, \Gamma_{\rm noise})$ is independent of $X$. Here, $\Gamma_{\rm pr} \in \R^{n \times n}$ and $\Gamma_{\rm noise} \in \R^{m \times m}$ are symmetric positive definite covariance matrices, $x_{\rm pr} \in \R^n$ is the prior mean for $X$, and we have explicitly indicated the nonlinear dependence of the discrete measurement matrix $R(p)$ on the positioning of the parallel beam source-receiver pair. Under these simplifying assumptions, the posterior in \eqref{eq:Bayes} is also Gaussian with the covariance matrix and mean \cite{kaipio2006statistical}
\begin{subequations}
   \label{eq:posterior_cov}
  \begin{align}
  \Gamma_{\rm post}(p) &= \big(\Gamma_{\rm pr}^{-1} + R(p)^T \Gamma_{\rm noise}^{-1} R(p) \big)^{-1},\\[1mm]
  \widehat{x}(p) &= \Gamma_{\rm post}(p) \big( \Gamma_{\rm pr}^{-1}x_{\rm pr} + R(p)^T \Gamma_{\rm noise}^{-1} y \big),
\end{align}
\end{subequations}
respectively, as can be deduced by a straightforward completion of squares in \eqref{eq:Bayes}. Using the Woodbury matrix identity, these equations can alternatively be written as
\begin{subequations}
   \label{eq:posterior_cov2}
  \begin{align}
  \Gamma_{\rm post}(p) &= \Gamma_{\rm pr} - \Gamma_{\rm pr} R(p)^T \big(R(p) \Gamma_{\rm pr} R(p)^T + \Gamma_{\rm noise}  \big)^{-1} R(p) \Gamma_{\rm pr} ,\\[1mm]
  \widehat{x}(p) &= x_{\rm pr} +  \Gamma_{\rm pr} R(p)^T \big(R(p) \Gamma_{\rm pr} R(p)^T + \Gamma_{\rm noise}  \big)^{-1} (y -  R(p) x_{\rm pr}),
  \end{align}
  which are computationally far more attractive than \eqref{eq:posterior_cov} if $m \ll n$, as is the case in our setting. However, if one is interested in the inverse covariance $\Gamma_{\rm post}(p)^{-1}$, then \eqref{eq:posterior_cov} should be used.
\end{subequations}

\subsection{A and D-optimality}
\label{sec:A_and_D}
In Bayesian optimal experimental design, one often considers minimizing the expected squared distance of the unknown in a given (semi)norm around some chosen point estimate, which corresponds to the so-called A-optimal design; in all our considerations, the point estimate of interest is the posterior mean. The other commonly used alternative is to look for the D-optimal design, that is, to maximize the information gain when the prior is replaced by the posterior. See,~e.g.,~\cite{alexanderian2016bayesian,chaloner1995bayesian} for more details.

Assuming  the employed seminorm is induced by the positive semidefinite matrix $A^T \!A$ for a given $A \in \R^{l \times n}$, in our simple, i.e.~Gaussian, linear and finite-dimensional, setting, A-optimality corresponds to choosing a design parameter $p_{\rm A} \in \R^2$ satisfying
\begin{equation}
\label{eq:Aoptimal}
p_{\rm A} = {\rm arg} \min_{p} \, {\rm tr}  \big(A \Gamma_{\rm post}(p) A^T\big).
\end{equation}
On the other hand, D-optimality is equivalent to finding
\begin{equation}
  \label{eq:Doptimal}
 p_{\rm D} = {\rm arg} \min_{p} \, \log \!\big(\!\det \Gamma_{\rm post}(p) \big),
\end{equation}
where the logarithm follows from the definition of D-optimality but it also makes the numerical treatment of \eqref{eq:Doptimal} considerably more stable as discussed in the following subsection. For completeness, the deductions of \eqref{eq:Aoptimal} and \eqref{eq:Doptimal} are presented in Appendix~\ref{sec:KL}.

Since all pixels are of the same size in our discretized model, choosing $A = I \in \R^{n\times n}$ to be the identity matrix leads to using (a scaled) $L^2(D)$-norm as the distance measure in \eqref{eq:Aoptimal}. If one is only interested in the expected $L^2$-accuracy of the posterior mean over some specific ROI within $D$, one simply only needs to replace the diagonal elements of $A$ corresponding to the pixels in the complement of the ROI by zeros. On the other hand, being only interested in the information gain in the ROI for D-optimality corresponds to replacing $\Gamma_{\rm post}(p)$ in \eqref{eq:Doptimal} by the matrix obtained by deleting the rows and columns of $\Gamma_{\rm post}(p)$ corresponding to the uninteresting pixels. We refer to Appendix~\ref{sec:KL} for the proof of this statement as well as for the exact form of the information gain when the prior is replaced by the posterior, as it equals  $- \log  \det (\Gamma_{\rm post}(p))$ only up to an affine transformation.

\subsection{Evaluation of the target functionals}
\label{sec:evaluation}
Since the target functionals in \eqref{eq:Aoptimal} and \eqref{eq:Doptimal} typically have multiple local minima, we choose a straightforward approach and perform an exhaustive search over a two-dimensional grid. In consequence, the ability to efficiently evaluate the optimization targets becomes the top priority. To this end, observe that the (effective) rank of the perturbation $R(p)^T \Gamma_{\rm noise}^{-1} R(p)$ in \eqref{eq:posterior_cov} is presumably $m$ since the measurement matrix $R(p)$ corresponding to a single parallel beam projection is expected to be of full rank with all its $m$ singular values being of the same order of magnitude (cf.~Section~\ref{sec:mm}). Hence, our plan is to employ the matrix determinant lemma and the alternative formula for the posterior covariance in \eqref{eq:posterior_cov} in order to only need to compute traces of inverses and log-determinants for matrices of size $m \times m$. In particular, there is no apparent reason for resorting to Monte Carlo techniques, such as those in \cite{alexanderian2016fast,alexanderian2018efficient,avron2011randomized}, because the low (effective) rank of the perturbation $R(p)^T \Gamma_{\rm noise}^{-1} R(p)$ in \eqref{eq:posterior_cov} can be exploited explicitly.

With a suitable numbering of the pixels, the posterior covariance matrix can be written in a block form as
\begin{equation}
\label{eq:block_form}
\Gamma_{\rm post} =
\begin{bmatrix}
	\Gamma_{\rm ROI} & \Gamma_{\rm mix}\\[0.5em]
	\Gamma_{\rm mix}^T & \Gamma_{\rm rest},
\end{bmatrix}
=
\begin{bmatrix}
	(\Gamma_{\rm post})_{11} & (\Gamma_{\rm post})_{12}\\[0.5em]
	(\Gamma_{\rm post})_{21} & (\Gamma_{\rm post})_{\rm 22}
\end{bmatrix} \in \R^{n\times n},
\end{equation}
where $\Gamma_{\rm ROI} = (\Gamma_{\rm post})_{11}$ and $\Gamma_{\rm rest} = (\Gamma_{\rm post})_{\rm 22}$ are the (marginal) covariance matrices for the pixels in the ROI and the rest of the image, respectively. We use a similar indexing for analogous block decompositions of other $n \times n$ matrices as well. Moreover, as in \eqref{eq:block_form}, we often simplify the notation by not explicitly marking the dependence on $p$.

\begin{remark}
Because our final product is a sequential optimization algorithm, the prior covariance for an optimization step is usually the optimized posterior from the previous one. Hence, having an explicit representation for $\Gamma_{\rm prior}^{-1}$ at hand is arguably the most natural assumption for a single step of the algorithm; see \eqref{eq:posterior_cov}. However, assuming $m$ is of moderate size, $\Gamma_{\rm prior}$ can also be formed explicitly via \eqref{eq:posterior_cov} without consuming too much time, and thus we assume to explicitly know $\Gamma_{\rm prior}$  in the following. On the other hand, if the noise covariance is diagonal, as it is in all our numerical experiments, knowing $\Gamma_{\rm noise}$ is essentially the same as knowing $\Gamma_{\rm noise}^{-1}$.
\end{remark}

\subsubsection{D-optimality}
As discussed in Section~\ref{sec:A_and_D}, the aim of D-optimal design with a preassigned ROI is to minimize
\begin{equation}
  \label{eq:Dtarget}
\Phi_{\rm D}(p) := \log \!\big(\!\det \Gamma_{\rm ROI}(p) \big)
\end{equation}
over $p$ in a certain subset of $\R^2$. Recall that our aim is to efficiently evaluate $\Phi_{\rm D}(p)$ at numerous $p \in \R^2$.

To begin with, observe that
\begin{align*}
\det(\Gamma_{\rm post}) &= \det
\begin{bmatrix}
	\Gamma_{\rm ROI} & \Gamma_{\rm mix}\\[0.5em]
	0 & \Gamma_{\rm rest} - \Gamma_{\rm mix}^T \Gamma_{\rm ROI}^{-1}\Gamma_{\rm mix}
\end{bmatrix}\\[1em]
	&= \det(\Gamma_{\rm ROI})\det(\Gamma_{\rm rest} - \Gamma_{\rm mix}^T \Gamma_{\rm ROI}^{-1}\Gamma_{\rm mix}),
\end{align*}
where $\Gamma_{\rm ROI}$ is invertible as a nonempty diagonal block of a positive definite matrix.
In particular, $\Gamma_{\rm rest} - \Gamma_{\rm mix}^T \Gamma_{\rm ROI}^{-1}\Gamma_{\rm mix}$
is the Schur complement of $\Gamma_{\rm ROI}$, meaning that
\[
\Gamma_{\rm post}^{-1} =: \Sigma =
\begin{bmatrix}
	\Sigma_{11} & \Sigma_{12}\\
	\Sigma_{21} & \Sigma_{22}
\end{bmatrix} =
\begin{bmatrix}
	\ * \ & \ * \ \\[1mm]
	\ * \  & (\Gamma_{\rm rest} - \Gamma_{\rm mix}^T \Gamma_{\rm ROI}^{-1}\Gamma_{\rm mix})^{-1}\\
\end{bmatrix}.
\]
Altogether we thus have
\[
	\det(\Gamma_{\rm ROI}) =
	\frac{\det(\Gamma_{\rm post})}{\det(\Gamma_{\rm rest} - \Gamma_{\rm mix}^T \Gamma_{\rm ROI}^{-1}\Gamma_{\rm mix})} =
	\frac{\det(\Sigma_{22})}{\det(\Sigma)},
\]
and so one needs to consider evaluating the logarithms of $\det(\Sigma)$ and $\det(\Sigma_{22})$, both of which depend on the design parameter $p$.

In the rest of this section, we explicitly mark which matrices depend on $p$. According to \eqref{eq:posterior_cov},
\begin{align*}
	\Sigma(p) &= \Gamma_{\rm pr}^{-1} + R(p)^{T} \Gamma_{\rm noise}^{-1}R(p),\\[1mm]
	\Sigma_{22}(p) &= (\Gamma_{\rm pr}^{-1})_{22} + R_2(p)^{T} \Gamma_{\rm noise}^{-1}R_2(p),
\end{align*}
where $R(p) = \left[ R_1(p) , R_2(p) \right] \in \R^{m \times n}$ is a columnwise partitioning, with $R_1(p)$ corresponding to the ROI and $R_2(p)$ to the rest of the image. By virtue of the matrix determinant lemma,
\begin{align*}
	\det(\Sigma(p)) &= \det\!\big( R(p) \Gamma_{\rm pr}R(p)^T + \Gamma_{\rm noise} \big) \det(\Gamma_{\rm noise}^{-1}  )\det(\Gamma_{\rm pr}^{-1}),\\[1mm]
	\det(\Sigma_{22}(p)) &= \det\!\big( R_2(p) \big((\Gamma_{\rm pr}^{-1})_{22}\big)^{-1} R_2(p)^T + \Gamma_{\rm noise} \big)
	\det(\Gamma_{\rm noise}^{-1}) \det\!\big((\Gamma_{\rm pr}^{-1})_{22}\big),
\end{align*}
where only the first determinants, respectively, depend on $p$. Moreover, the matrices associated to these $p$-dependent determinants are only of size $m \times m$,~i.e.~small.

To evaluate $\log \det(\Gamma_{\rm ROI}(p))$, one thus needs to precompute Cholesky decomposition for the (large) positive definite matrix $(\Gamma_{\rm pr}^{-1})_{22}$ that is independent of $p$. Since $m$, i.e.~the number of columns in $R(p)^T$, is assumed to be relatively low, this further enables building the (small) Cholesky decompositions
\begin{subequations}
    \label{eq:chol}
\begin{align}
  C(p) C(p)^T &= R(p) \Gamma_{\rm pr}R(p)^T + \Gamma_{\rm noise} , \\
\tilde{C}(p) \tilde{C}(p)^T &= R_2(p) \big((\Gamma_{\rm pr}^{-1})_{22}\big)^{-1} R_2(p)^T  + \Gamma_{\rm noise}
\end{align}
\end{subequations}
for all $p$ on the employed grid. Finally, a straightforward algebraic manipulation gives
\begin{equation}
  \label{eq:dopt_final}
\log\!\big(\!\det(\Gamma_{\rm ROI}(p)) \big)=  2 \sum_{j=1}^m \big( \log (\tilde{c}_{jj}(p)) - \log (c_{jj}(p)) \big) + c,
\end{equation}
where $c_{jj}(p)$ and $\tilde{c}_{jj}(p)$, $j=1, \dots, m$, are the diagonal elements of $C(p)$ and $\tilde{C}(p)$, respectively, and $c \in \R$ is independent of $p$.

Observe that it is far more stable to numerically evaluate the sum of logarithms in \eqref{eq:dopt_final} than the products of the diagonal elements of the Cholesky factors needed for computing $\det(\Gamma_{\rm ROI}(p))$ itself. Moreover, the constant $c$ in \eqref{eq:dopt_final} can be evaluated by considering Cholesky decompositions for $(\Gamma_{\rm pr}^{-1})_{22}$ and  either $\Gamma_{\rm pr}$ or $\Gamma_{\rm pr}^{-1}$; the former was already employed when building \eqref{eq:chol}, and the latter neither poses any difficulties if $n$ is not huge. In fact, one can even evaluate the actual information gain when the prior is replaced by the posterior for any $p$ without too severely compromising the computational efficiency; see Appendix~\ref{sec:KL} for further details. This makes it possible to compare the information gains between different rounds of the sequential optimization algorithm introduced in Section~\ref{sec:opti} below.

\subsubsection{A-optimality}
As discussed in Section~\ref{sec:A_and_D}, when aiming at an A-optimal design, one needs to minimize
$$
\Phi_{\rm A}(p) := {\rm tr} (A \Gamma_{\rm post}(p) A^T)
$$
over $p$ in a certain subset of $\R^2$ for a given $A \in \R^{l \times n}$. In our numerical tests presented in Section~\ref{sec:numerics}, $A = I_{\rm ROI } \in \R^{n \times n}$ is an identity matrix with the diagonal elements corresponding to the complement of the ROI replaced by zeros. Nevertheless, we perform the following calculations for an unspecified $A$ to maintain generality.

As in the case of D-optimality, we introduce a Cholesky decomposition
\begin{equation*}
C(p) C(p)^T = R(p) \Gamma_{\rm pr}R(p)^T + \Gamma_{\rm noise} \in \R^{m \times m},
\end{equation*}
and then form an auxiliary matrix
$$
B(p) = C(p)^{-1} R(p) \Gamma_{\rm pr} A^T \in \R^{m \times l}.
$$
According to \eqref{eq:posterior_cov2},
$$
\Phi_{\rm A}(p) = {\tr} \big( A \Gamma_{\rm pr}A^T - B(p)^T \! B(p) \big)  = {\tr} ( A \Gamma_{\rm pr}A^T) - {\tr} (B(p)^T \! B(p)).
$$
Hence, fundamental properties of the matrix trace allow the representation
$$
\Phi_{\rm A}(p) = c' - \sum_{j=1}^m \sum_{k=1}^{l} B_{jk}(p)^2 .
$$
Since $c' = {\rm tr}(A \Gamma_{\rm pr} A^T)$ does not depend on $p$, it does not affect the minimization of $\Phi_{\rm A}(p)$. However, evaluating $c'$ does not considerably slow down the minimization process as a whole because it needs only be done once; knowing $c'$ allows comparing the minimal values of $\Phi_{\rm A}$ between different iterations of the sequential minimization algorithm introduced in the following section.

\section{Sequential optimization of measurements}
\label{sec:opti}
In this section we present the algorithm for sequentially optimizing the parallel beam projections for X-ray tomography. We start with the basic algorithm that aims at finding the optimal projections prior to having any measurements at hand. Subsequently, we consider the modifications needed if one wants to estimate  a hyperparameter in the prior based on the observed data or to change the ROI for the $(k+1)$th optimization round based on the reconstruction, i.e.~the {\em conditional mean} (CM) estimate, computed from the first $k$ projections.

\subsection{Basic algorithm}
The algorithm is initialized by determining the discretization of $D$ into $n$ pixels, choosing the spacing of the $m$ detectors/rays for a single parallel beam projection, and specifying the covariance matrices $\Gamma_{0}$ and $\Gamma_{\rm noise}$ for the Gaussian prior of the absorption coefficient and for the Gaussian noise model, respectively.\footnote{The mean of the prior is not needed for the sequential optimization algorithm, unless it affects the noise model.} One also needs to choose the ROI, select the number of optimization rounds, and determine how the design variable $p \in \R^2$  parametrizes a parallel beam projection.  In our numerical experiments, the first component of $p$ corresponds to the projection angle and the second one to the distance of the bisection of the parallel beam source-receiver pair from the origin.

The optimization is performed sequentially, i.e., so that the posterior probability density after the $k$th projection is used as the prior for choosing the $(k+1)$th projection. At each step the target functional for A or D-optimality is evaluated for all $p$ on a two-dimensional optimization grid that also needs to be predetermined. The particular $p$ yielding the minimum value for the considered minimization target is chosen as the optimal parameter, and subsequently the posterior covariance, i.e.~the prior covariance for the next optimization step, is formed according to \eqref{eq:posterior_cov} or \eqref{eq:posterior_cov2}.

The algorithm proceeds altogether as follows:
\begin{algorithm}
\label{alg:basic_optimization}
\begin{algorithmic}
  \STATE{Choose the covariances $\Gamma_0$ and $\Gamma_{\rm noise}$ for the initial Gaussian prior $\mathcal{N}(x_0, \Gamma_0)$ and the noise model $\mathcal{N}(0, \Gamma_{\rm noise})$. Select the ROI, the grid for $p$ and the number of optimization rounds $K$. Set $\Gamma_{\rm pr} = \Gamma_0$ and $P = [ \ \ ]$.}
  \FOR{$k=1, \dots, K$}
    \FOR{all $p$ on the optimization grid}
      \STATE{Evaluate $\Phi_{\rm D}(p)$ or $\Phi_{\rm A}(p)$ as outlined in Section~\ref{sec:evaluation}.}
    \ENDFOR
    \STATE{Find the minimizer $p_k$ of the considered optimization target.}
    \STATE{Append $P = [P, \ p_k]$.}
    \STATE{Form the posterior covariance $\Gamma_{\rm post}(p_k)$.}
    \STATE{Set $\Gamma_{\rm pr} = \Gamma_{\rm post}(p_k)$.}
  \ENDFOR
  \RETURN $P$
\end{algorithmic}
\end{algorithm}

The columns of the output matrix $P \in \R^{2 \times K}$ define the $K$ sequentially optimized parallel beam projections. To be more precise, the $(k+1)$th column of $P$ defines the A or D-optimal projection given the previous $k$ projections, that is, each of the optimized projections is (only) {\em locally} optimal. In particular, there is no reason to expect that the found $K$ parallel beam projections would be {\em globally optimal} \cite{huan2016sequential}: It is a simpler and computationally less demanding task to optimize the projections one by one than to simultaneously find $K$ parallel beam projections that are jointly A or D-optimal. However, there {\em is} anyway reason to expect that the projections defined by the columns of $P$ are more optimal than, e.g., a random choice.

Due to the assumptions that the measurement model is linear with additive noise and the prior and noise process are independent Gaussians, Algorithm \ref{alg:basic_optimization} can be run prior to performing any measurements. Hence, one may expect to have lots of computational time and resources for performing the sequential optimization of Algorithm \ref{alg:basic_optimization} in an `offline mode'. However, it is also possible to change the ROI for the next optimization step or to estimate some free parameters in the prior based on the previous data or reconstruction as explained in the following two subsections.

\subsection{Adaptive region of interest}
If the sequentially optimal angles are not determined before the measurements but as a part of the online imaging procedure, the ROI may be altered between the optimization rounds based on the observations of the expert running the algorithm and,~e.g.,~interesting or alarming features in the previous reconstruction. In this case, one also needs to give the initial prior mean as an input and incorporate the computation of the CM estimate in each step of the algorithm; see \eqref{eq:posterior_cov} and \eqref{eq:posterior_cov2}. Moreover, the natural output is no longer the sequence of optimal projections but the final reconstruction, and the outer iterations should be stopped by the operator of the algorithm.

\begin{algorithm}[adaptive ROI]
\label{alg:ROI_optimization}
\begin{algorithmic}
  \STATE{Choose the covariances $\Gamma_0$ and $\Gamma_{\rm noise}$ as well as the mean $x_0$ for the initial Gaussian prior $\mathcal{N}(x_0, \Gamma_0)$ and the noise model $\mathcal{N}(0, \Gamma_{\rm noise})$. Select the ROI and the optimization grid for $p$. Initialize $\Gamma_{\rm pr} = \Gamma_0$, $x_{\rm pr} = x_0$ and $k=0$.}
  \WHILE{satisfactory reconstruction has not been reached}
  \STATE{Set $k \leftarrow k+1$.}
    \FOR{all $p$ on the optimization grid}
      \STATE{Evaluate $\tr(I_{\rm ROI} \Gamma_{\rm post}(p) I_{\rm ROI})$ or $\log( \det(\Gamma_{\rm ROI}(p)))$ as outlined in Section~\ref{sec:evaluation}.}
    \ENDFOR
    \STATE{Find the minimizer $p_k$ of the considered optimization target.}
    \STATE{Observe new data $y = R(p_k) x + n_k$, where $n_k$ is a realization of the noise.}
    \STATE{Form the posterior covariance $\Gamma_{\rm post}(p_k)$ and the CM estimate $\widehat{x}(p_{k})$.}
    \STATE{Redefine (heuristically) the ROI based on $\widehat{x}(p_k)$.}
    \STATE{Set $\Gamma_{\rm pr} = \Gamma_{\rm post}(p_k)$ and $x_{\rm pr} = \widehat{x}(p_{k})$.}
    \ENDWHILE
    \STATE{Set $x_{\rm rec} = x_{\rm pr}$ and $\Gamma_{\rm rec} = \Gamma_{\rm pr}$.}
  \RETURN $x_{\rm rec}$ and $\Gamma_{\rm rec}$
\end{algorithmic}
\end{algorithm}

It follows from basic Bayesian analysis that the sequentially updated reconstruction $x_{\rm rec} \in \R^n$ and the spread estimator $\Gamma_{\rm rec}$ produced by Algorithm~\ref{alg:ROI_optimization} are the same one would obtain by first collecting the data for all optimized projection angles and only then computing the CM estimate and the posterior covariance in a single step.

\subsection{Data-driven hyperparameter estimation}
Suppose the {\em original} prior covariance $\Gamma_{0} = \Gamma_{0}(\rho) \in \R^{n \times n}$ depends on a hierarchical parameter $\rho \in \R^l$, the true value of which is unknown. Our aim is to introduce an algorithm that alternates between determining the optimal design parameter and the maximum likelihood estimate for the hyperparameter based on the measured data. The latter step can be considered as empirically adapting the prior at the same time as the sequential design process is performed.

Denote by
\begin{equation*}
	\Gamma_{\rm ms}(\rho;  p) = \Gamma_{\rm noise} + R(p) \Gamma_{0}(\rho) R(p)^T \in \R^{m \times m}
\end{equation*}
the covariance matrix of the measurement given the prior covariance $\Gamma_{0}(\rho)$ and the projection matrix $R(p)$. The marginal density of the measurement $y$ conditioned on $\rho$ is hence given by
\begin{equation*}
	\pi(y \, | \, \rho;  p) \propto
	\exp \!\Big(\! - \! \frac{1}{2} \big(  \log (\det \Gamma_{\rm ms}(\rho ; p)) - (y - R(p) x_0 )^T   \Gamma_{\rm ms}(\rho; p)^{-1}  (y - R(p) x_0 ) \big) \Big),
\end{equation*}
where $x_0 \in \R^n$ is the (original) prior mean. In consequence, having $k \in \N$ independent measurements $y_1, \dots, y_k$ corresponding to the design parameters $p_1, \dots, p_k$ at hand leads to
\begin{align}
  \label{eq:likelihood}
	\pi&(y_1, \dots , y_k \, | \, \rho;  p_1, \dots, p_k) = \prod_{i=1}^k \pi(y_i \, | \, \rho;  p_i) \\[-1mm]
& \propto \exp\! \Big(\! - \! \frac{1}{2} \sum_{i=1}^k \big(  \log (\det \Gamma_{\rm ms}(\rho ; p_i)) - (y_i - R(p_i) x_0 )^T   \Gamma_{\rm ms}(\rho; p_i)^{-1}  (y_i - R(p_i) x_0 ) \big) \Big) \nonumber.
\end{align}
After obtaining the latest measurement $y_k$ corresponding to the design parameter~$p_k$, one can thus determine the latest {\em maximum likelihood} (ML) estimate $\rho_k$ for the hyperparameter by maximizing \eqref{eq:likelihood} with respect to $\rho$.

In addition to the computational cost corresponding to maximizing \eqref{eq:likelihood}, an extra price one needs to pay for the empirical estimation of a hyperparameter is the need to appropriately update the posterior for $X$ to make it compatible with the current ML estimate $\rho_k$. In other words, one needs to recompute the posterior from a scratch:  $\Gamma_{\rm post}(\rho_k; p_1, \dots, p_k)$ is formed,~e.g.,~via \eqref{eq:posterior_cov} or \eqref{eq:posterior_cov2} with $\Gamma_{\rm pr} = \Gamma_0(\rho_k)$, $\Gamma_{\rm noise}$ built from $k$ diagonal blocks, and $R(p)$ replaced by the `total projection matrix' $[R(p_1)^T, \dots, R(p_k)^T]^T$.

\begin{algorithm}[data-driven hyperparameter estimation]
\label{alg:data-driven_adaptive_estimation}
\begin{algorithmic}
  \STATE{Choose the covariance $\Gamma_{\rm noise}$ for the noise model $\mathcal{N}(0, \Gamma_{\rm noise})$ and the mean $x_0$ for the Gaussian prior. Initialize $\rho_0 \in \R^l$, $\Gamma_{\rm pr} = \Gamma_0(\rho_0)$, and set $k=0$. Select the ROI and the grid for $p$.}
  \WHILE{e.g., a satisfactory reconstruction has not been reached}
  \STATE{Set $k \leftarrow k+1$.}
    \FOR{all $p$ on the optimization grid}
      \STATE{Evaluate $\Phi_{\rm D}(p)$ or $\Phi_{\rm A}(p)$ as outlined in Section~\ref{sec:evaluation}.}
    \ENDFOR
    \STATE{Find the minimizer $p_k$ of the considered optimization target.}
    \STATE{Observe the data $y_k = R(p_k) x + n_k$, where $n_k$ is a realization of the noise.}
  \STATE{Solve for the \rm ML estimate $\rho_k = {\rm arg\,max}_{\rho \in \R^l} \pi(y_1, \dots, y_k \, | \, \rho, p_1, \dots, p_k).$}
    \STATE{Form the posterior covariance $\Gamma_{\rm post}(\rho_k; p_1, \dots, p_k)$ ``from a scratch''.}
    \STATE{Set $\Gamma_{\rm pr} = \Gamma_{\rm post}(\rho_k; p_1, \dots, p_k)$.}
  \ENDWHILE
\end{algorithmic}
\end{algorithm}

Computing the ML estimate $\rho_k$ is relatively straightforward: When $k=1$, one needs to pay some attention to make sure the employed minimization technique for finding $\rho_1$ is not very slow and does not diverge in case there is not much information on the whereabouts of the true parameter value to begin with. However, when $\rho_k$ for $k \in \N$ is available, it can be used as the initial guess for finding $\rho_{k+1}$ by the Newton's method that exhibits fast convergence in our numerical experiments, where $\rho$ is the correlation length for a certain parametrized covariance structure. If the elementwise derivatives of $\Gamma_{0}(\rho)$ are known, the required first two derivatives for the term in the exponent of $\pi(y_1, \dots, y_k \, | \, \rho, p_1, \dots, p_k)$ can be straightforwardly evaluated by resorting to differentiation formulas for the determinant and the matrix inverse. Altogether, finding the needed ML estimates poses no severe computational difficulties if $m$ is relatively small and $l=1$, as is the case in our numerical studies.

\section{Numerical experiments}
\label{sec:numerics}
In all our tests, $D = [0,1]^2$ is the unit square divided into a uniform mesh of $n = N^2$ pixels. The maximal distance of X-rays passing through $D$ from its center point $(0.5,0.5)$  is set to $0.5$. If $D$ contains an obstruction that blocks X-rays, the rows corresponding to the lines passing through the obstruction are deleted from the considered projection matrices,  as the corresponding measurements do not carry any information on the unknown absorption and can thus be ignored. Hence, the number of active detectors may vary as a function of the position of the parallel beam source-receiver pair relative to the possible obstruction, which is taken into account when optimizing the optimal projections and computing related reconstructions. (It should be intuitively clear that optimal parallel beam projections try to avoid imaging through an obstruction whenever possible.)

The (initial) prior covariance between the remaining pixels is assumed to be of the form
\begin{equation}
 \label{eq:prior_cov}
 (\Gamma_{0})_{i,j} = \gamma^2 \exp \left(-\frac{| x_i - x_j |^2}{2\ell^2} \right),
 \end{equation}
where $| \cdot |$ denotes the Euclidean norm, $\ell>0$ is the so-called correlation length and $\gamma>0$ is the pixelwise standard deviation.  Here $x_i$ and $x_j$ are the coordinates of the pixels (or more specifically the pixel centers) with indices $i$ and $j$. Loosely speaking, the larger $\ell$ is, the less oscillating are random draws from the (initial) prior $\mathcal{N}(x_0, \Gamma_0)$. The components of the additive noise in \eqref{eq:meas_model} are assumed to be independent with a common standard deviation $\sigma > 0$.

\subsection{Test~1:~Basic algorithm}
Our first numerical tests consider the basic version of our optimization routine,~i.e.~Algorithm~\ref{alg:basic_optimization}. The aim is to demonstrate that the algorithm generates sequentially A-optimal parallel beam projections that are intuitively acceptable and clearly outperform random selections. Moreover, according to our tests, the optimization of projections can be performed using a considerably coarser discretization compared to the one used for computing the actual reconstructions, without significantly compromising the overall performance of the approach. Regarding D-optimality, it is numerically demonstrated that the sequential optimization scheme may sometimes lead to globally suboptimal sets of projections, although in many  cases the sets of sequentially A and D-optimal projections are actually qualitatively very similar. For completeness it should also be admitted that the possibility that sequentially A-optimal projections may also sometimes suffer from apparent global nonoptimality cannot be excluded based solely on the presented results.

\subsubsection*{Test 1.1: Sanity check}
\begin{figure}
  \begin{center}
		\hspace{-1cm}
	  \includegraphics[width = 0.53\columnwidth]{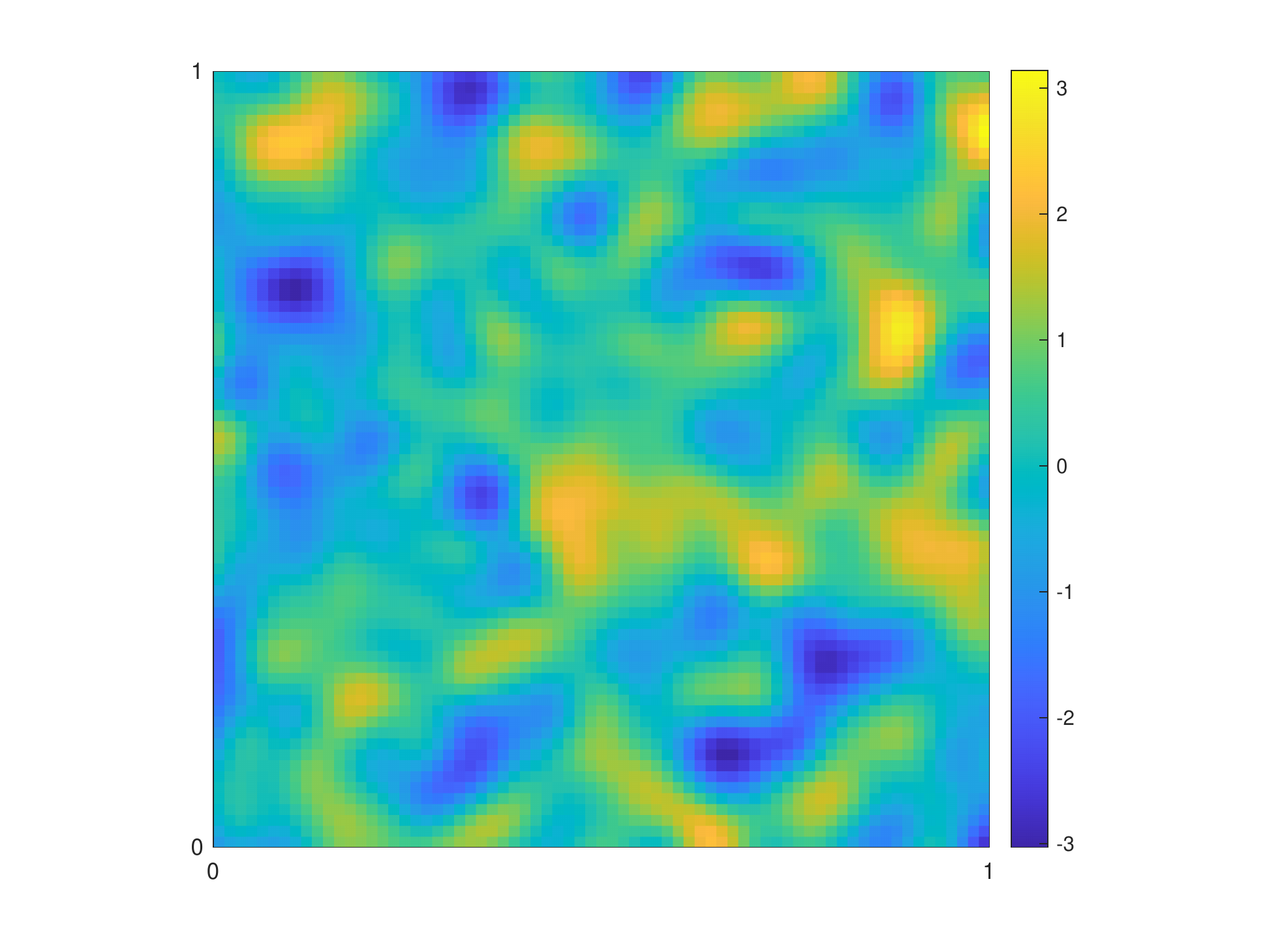}
		\hspace{-0.5cm}
	  \includegraphics[width = 0.53\columnwidth]{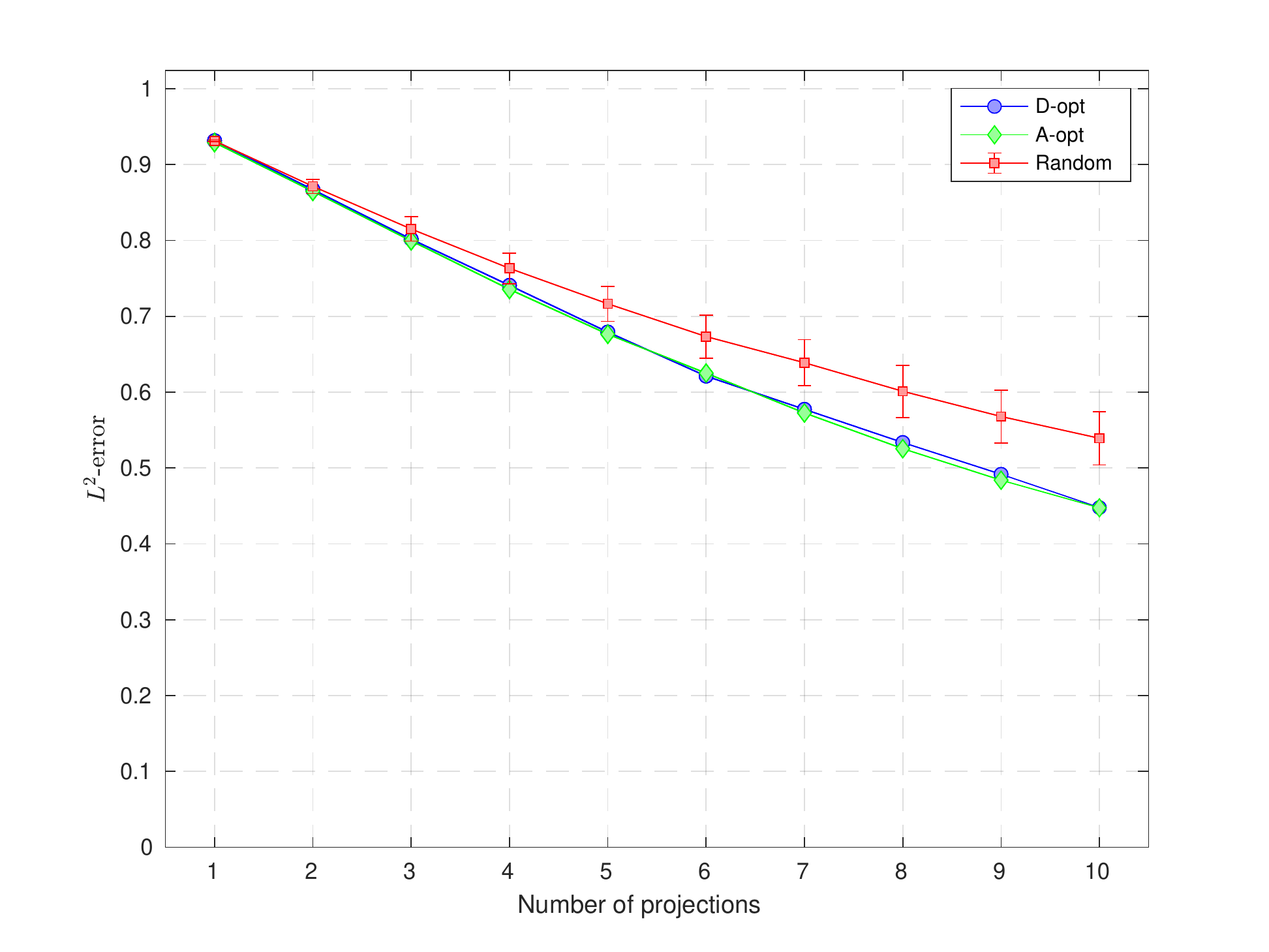}
  \end{center}
  \caption{{\sc Test~1.1.} Left: Random draw from the zero-mean Gaussian prior with the covariance structure \eqref{eq:prior_cov} for $\gamma = 1$ and $\ell = 0.05$.  Right: Mean $L^2(D)$ reconstruction errors over a sample of 1000 absorption targets as functions of the number of (noisy) projections. The diamonds correspond to the A-optimal projections and the circles to the D-optimal ones.
    The errorbars give the averages and standard deviations of the corresponding mean $L^2(D)$ reconstructions errors for 1000 sequences of projections with random angles.}
      \label{fig:test1_1b}
  \end{figure}

Let us start by applying Algorithm~\ref{alg:basic_optimization} to the simplest possible setting: the ROI is $D$, the parallel beam source-receiver pair is of the maximal width $1$, and there are no obstructions inside $D$. The number of pixels per edge of $D$ is chosen to be $N=100$ corresponding to altogether $n = 10^4$ pixels, and the number of individual detectors is $m=45$ per a single projection image. The pointwise standard deviation of the prior covariance in \eqref{eq:prior_cov} is assumed to be $\gamma = 1$ and the correlation length is chosen as $\ell = 0.05$. A random draw from the prior probability distribution of the absorption is visualized in the left-hand imaged of Figure~\ref{fig:test1_1b}. Here and in what follows, we assume the prior has zero mean, which is physically unrealistic as such, but mathematically this assumption simply corresponds to translating the coordinate system in $\R^n$ by some given physically sensible, component-wise positive prior mean $x_0$. The standard deviation of the additive white noise is $\sigma = 0.05$.

\begin{figure}
  \begin{center}
  \includegraphics[width = 0.48\columnwidth]{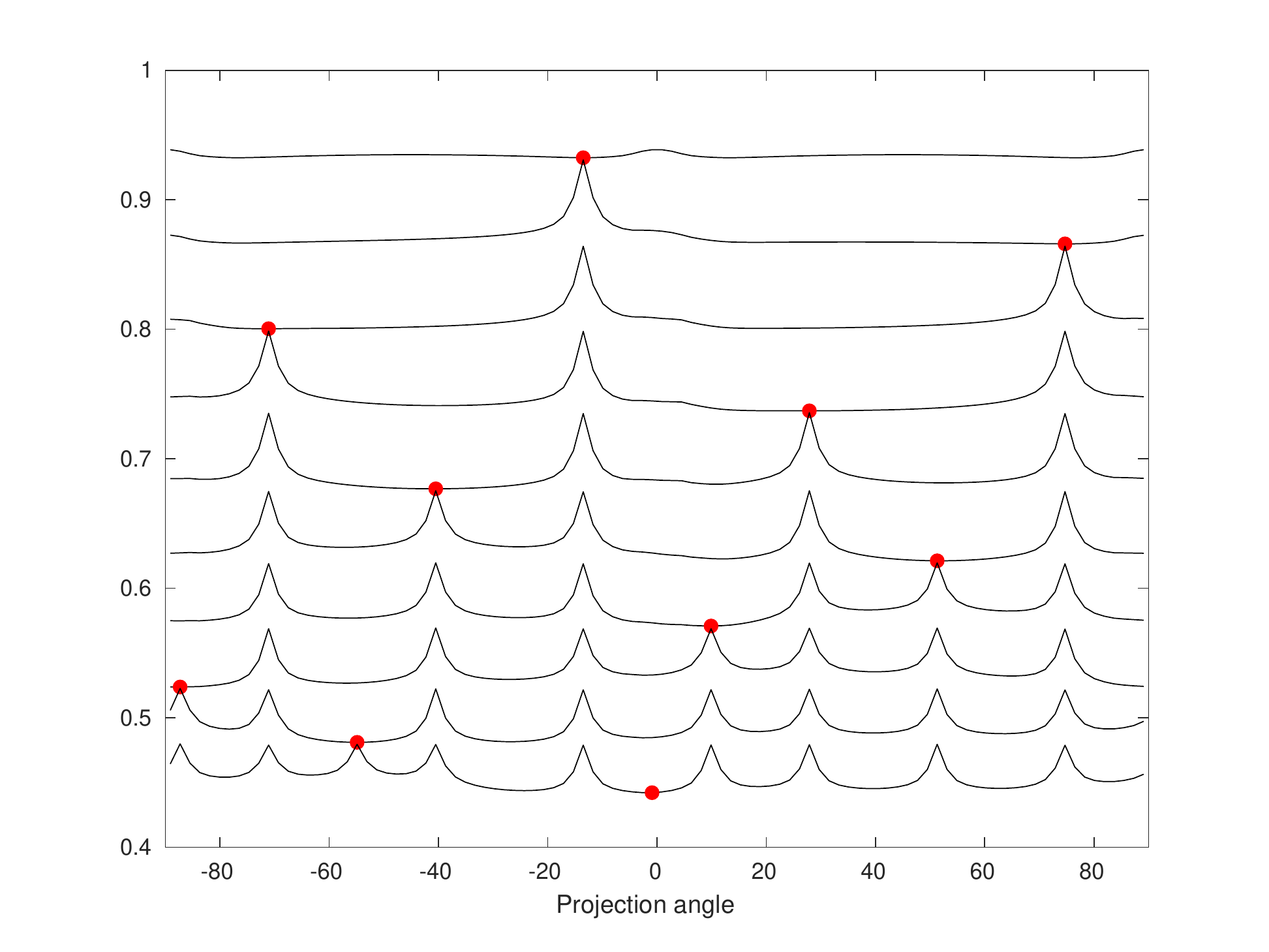}
  \includegraphics[width = 0.48\columnwidth]{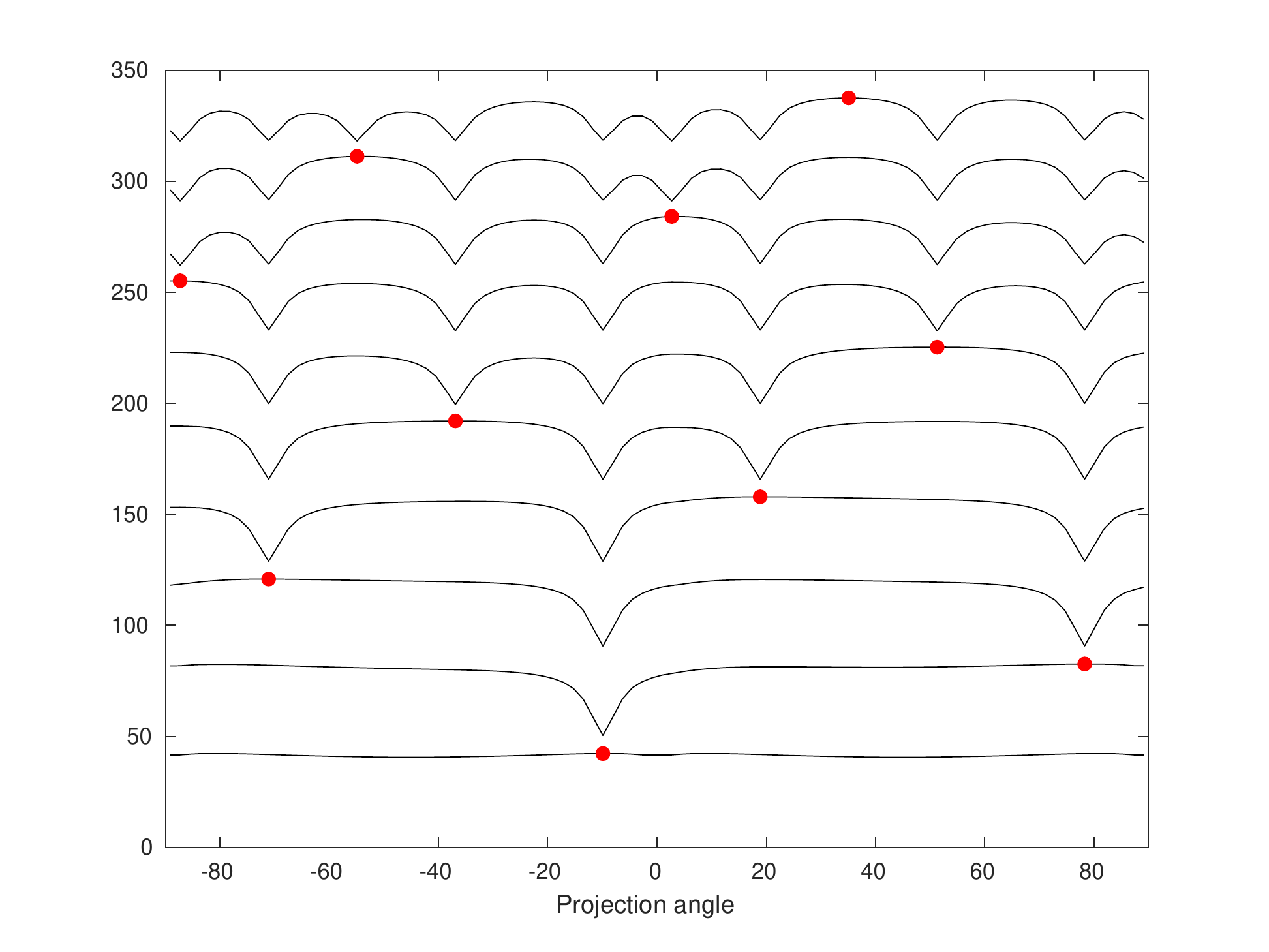}
  \end{center}
  \caption{{\sc Test~1.1.} Left: Modified A-optimality target $\frac{1}{N}\sqrt{\Phi_{\rm A}}$, corresponding to the expected $L^2(D)$ reconstruction error, as a function of the projection angle for the first ten projections. The dots denote the locally A-optimal angles and the algorithm proceeds from top to bottom (cf.~Figure~\ref{fig:test1_1b}). Right: Information gain $-\tilde{\Phi}_{\rm D}$ (cf.~\eqref{eq:Dtarget2}) as a function of the projection angle for the first ten projections. The dots denote the optimal angles and the algorithm proceeds from bottom to top.}
   \label{fig:test1_1a}
\end{figure}

Figure~\ref{fig:test1_1a} presents the A and D-optimality target functionals as functions of the projection angle over the first ten iterations of Algorithm~\ref{alg:basic_optimization}. To be quite precise, a modified A-optimality target $\frac{1}{N}\sqrt{\Phi_{\rm A}}$ is considered as it corresponds to the expected $L^2(D)$ reconstruction error, and in case of D-optimality, the depicted quantity is the actual information gain when the prior is replaced by the posterior. (The maximization of the information gain is equivalent to minimizing $\Phi_D$, as explained in Appendix~\ref{sec:KL}.) For both optimality criteria, the optimal angles are distributed rather uniformly over the interval $[-90,90]^\circ$ in such a way that the most recent projection angle is always located around the midpoint of the widest angular interval with no previous projections, which is intuitively what one would expect. It is noteworthy that after two projections the target functionals have multiple local optima; there is no reason to expect this would not be the case for more complicated imaging configurations as well.

The right-hand image of Figure~\ref{fig:test1_1b} shows the mean $L^2(D)$ reconstruction errors after each of the first ten A and D-optimal projections over a sample of $1000$ target absorptions drawn from the assumed prior. To be more precise, the noisy measurements are first simulated for each of the target absorptions and all optimal angles. Then, the reconstructions, i.e.~posterior means, are formed and compared to the corresponding targets when including the simulated noisy projections in the reconstruction process one by one in the same order as they were introduced by Algorithm~\ref{alg:basic_optimization}.  The A and D-optimal angles seem to perform equally well on average, even though the mean $L^2(D)$ reconstruction error is precisely the quantity the A-optimal angles were designed to minimize.

For comparison, random projection angles are also considered: the mean $L^2(D)$ reconstruction errors over the considered sample of absorptions are computed for 1000 sequences of random projection angles, the components of which are picked independently from the uniform distribution over $[-90,90]^\circ$. The right-hand image of Figure~\ref{fig:test1_1b} illustrates the averages of the resulting mean $L^2(D)$ reconstruction errors together with their standard deviations over the considered random sequences of projection angles. According to Figure~\ref{fig:test1_1b}, the A and D-optimal projections clearly outperform the random selections --- at least in this simple geometric setting.

\subsubsection*{Test~1.2:~Discoidal ROI}
\begin{figure}
  \begin{center}
		\input{figures_final/circle}
	  \qquad
	  \includegraphics[width = 0.45\columnwidth]{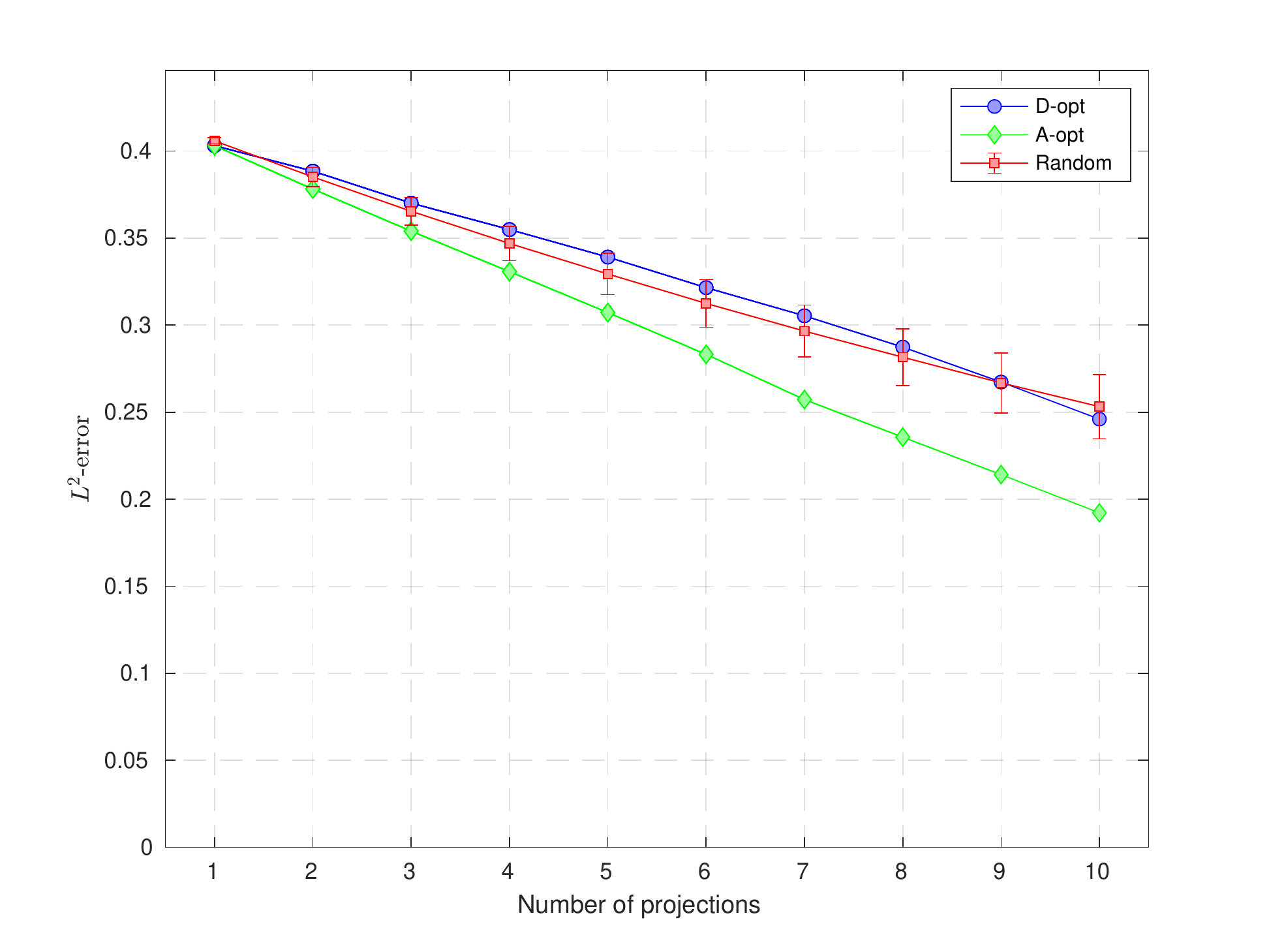}
  \end{center}
  \caption{{\sc Test~1.2.} Left: ROI. Right: Mean $L^2(D)$ reconstruction errors over a sample of 1000 absorption targets as functions of the number of (noisy) projections. The diamonds correspond to the A-optimal projections and the circles to the D-optimal ones. The errorbars give the averages and standard deviations of the corresponding mean $L^2(D)$ reconstructions errors for 1000 sequences of random projection that cover the ROI.}
      \label{fig:test1_2c}
\end{figure}

The second test is essentially a repetition of the first one, but with a parallel beam source-receiver pair of narrower width $0.5$ and the ROI specified to be a disk of radius $0.25$ centered at $(0.6,0.6)$; see the left-hand image of Figure~\ref{fig:test1_2c}. Observe, in particular, that the positioning of the source-receiver pair can now also be adjusted laterally and its width is the same as the diameter of the ROI. The standard deviation of the noise is set to $\sigma = 0.02$, and the number of individual detectors in the narrower parallel beam sensor is $m=23$. The other parameters are as in Test 1.1, that is, $N=100$, $\gamma = 1$ and $\ell = 0.05$.

\begin{figure}
  \begin{center}
		\includegraphics[width=\columnwidth]{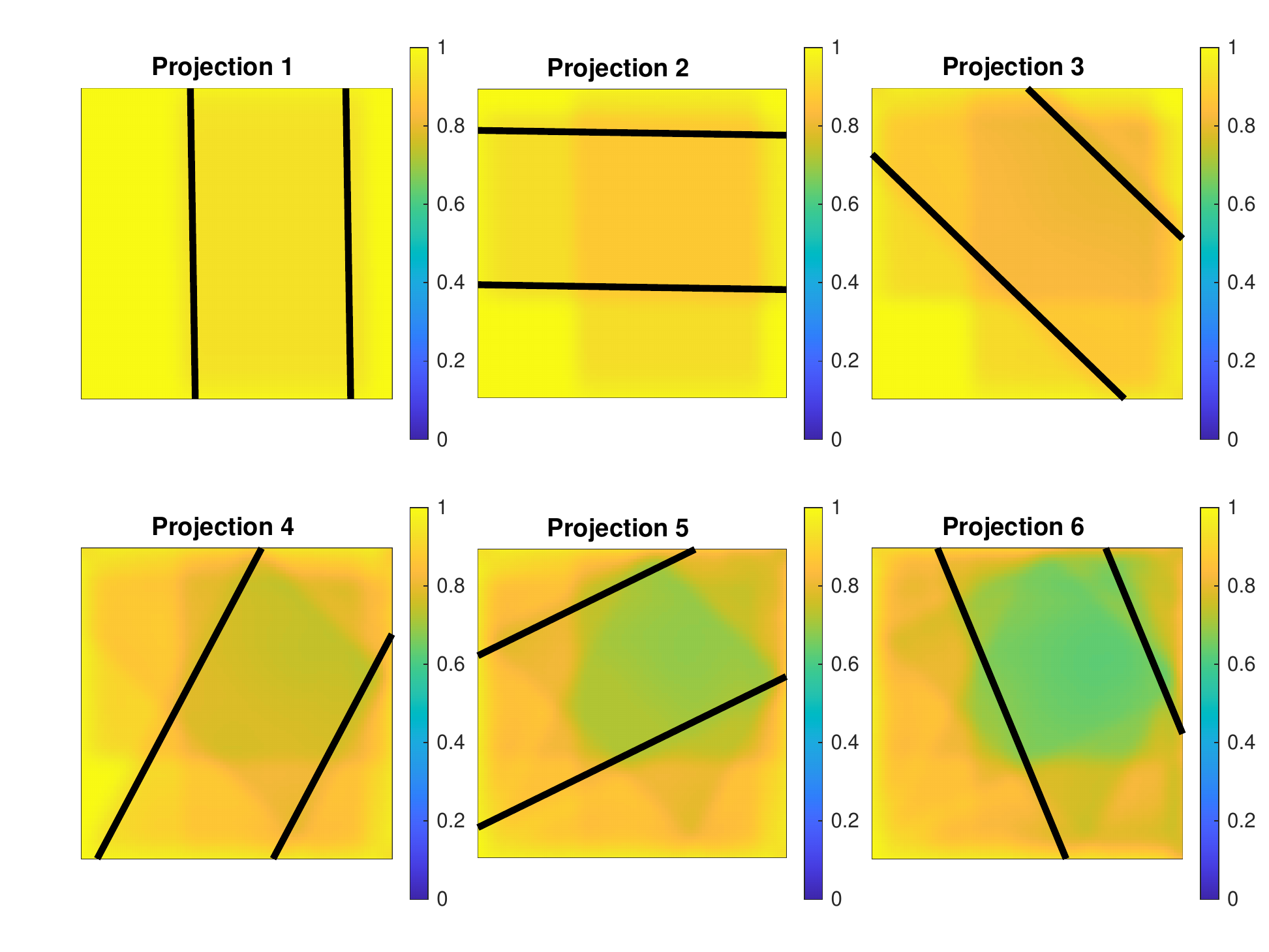}
  \end{center}
  \caption{{\sc Test~1.2.} The first six sequentially A-optimal parallel beam projections with the corresponding pixelwise (posterior) standard deviations in the background.}
      \label{fig:test1_2a}
  \end{figure}

\begin{figure}
  \begin{center}
		\includegraphics[width=\columnwidth]{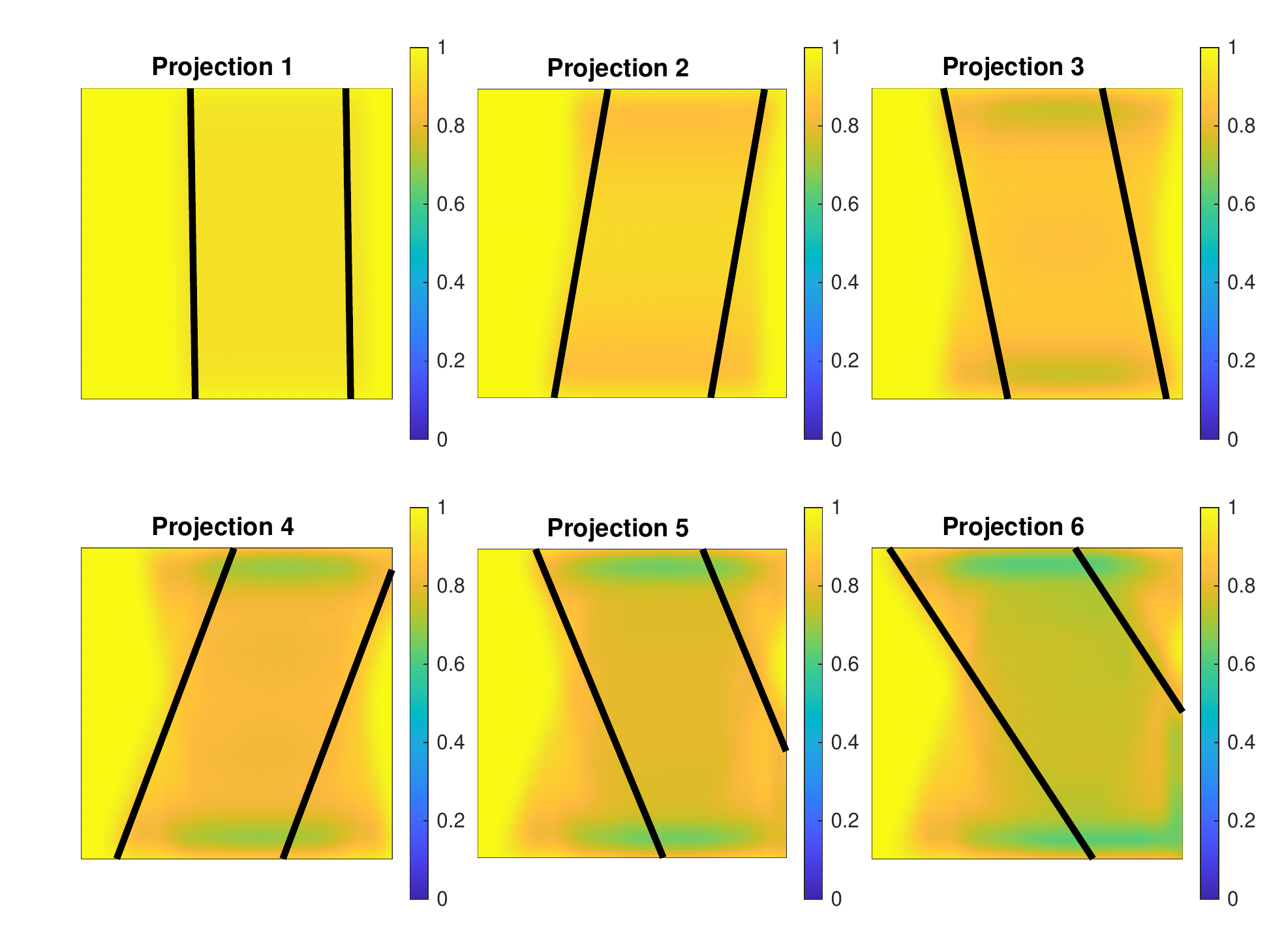}
  \end{center}
  \caption{{\sc Test~1.2.} The first six sequentially D-optimal parallel beam projections with the corresponding pixelwise (posterior) standard deviations in the background.}
      \label{fig:test1_2b}
\end{figure}

Figures~\ref{fig:test1_2a} and \ref{fig:test1_2b} visualize the first six sequentially A and D-optimal projections, respectively, with the corresponding pixelwise (posterior) standard deviations shown in the background. The A-optimal projections behave as intuitively as in Test~1.1: all projections approximately cover the ROI and they are also distributed relatively uniformly over all angles. Although the D-optimal projections also cover the ROI, they correspond to subsequent rotations of about $10^\circ$ relative to one of the previous projections. Even though such projections are indeed locally D-optimal, they definitely do not form a globally D-optimal set of projections: for example, the six sequentially A-optimal projections in Figure~\ref{fig:test1_2a} jointly produce a lower value for the D-optimality target $\Phi_{\rm D}$ than the sequentially D-optimal ones in Figure~\ref{fig:test1_2b}. This demonstrates that our sequential optimization procedure does not always produce parallel beam projections that are close to the global optimum.

The right-hand image of Figure~\ref{fig:test1_2c} illustrates the mean $L^2(D)$ reconstruction errors after each of the first ten (noisy) A and D-optimal projections over a sample of $1000$ target absorptions drawn from the assumed prior. As in Figure~\ref{fig:test1_1b} of Test~1.1, these mean errors are compared with the averages and standard deviations of the corresponding mean $L^2(D)$ reconstruction errors for $1000$ (semi)randomly selected projection sequences. To be more precise, the angles of the projections in these sequences are picked randomly from the uniform distribution over $[-90, 90]^\circ$, but the lateral position of the parallel beam source-receiver pair is subsequently adjusted so that all projections cover the ROI in order to allow a fair comparison with the optimized projections. The sequentially A-optimal projections clearly produce the lowest mean $L^2(D)$ reconstruction errors, but this time around the random selections slightly outperform the sequentially D-optimal projections in the considered performance metric.

\subsubsection*{Test~1.3: Obstruction and effect of coarse discretization}
\begin{figure}
  \begin{center}
		\input{figures_final/one_obstruction}
	  \qquad
	  \includegraphics[width = 0.45\columnwidth]{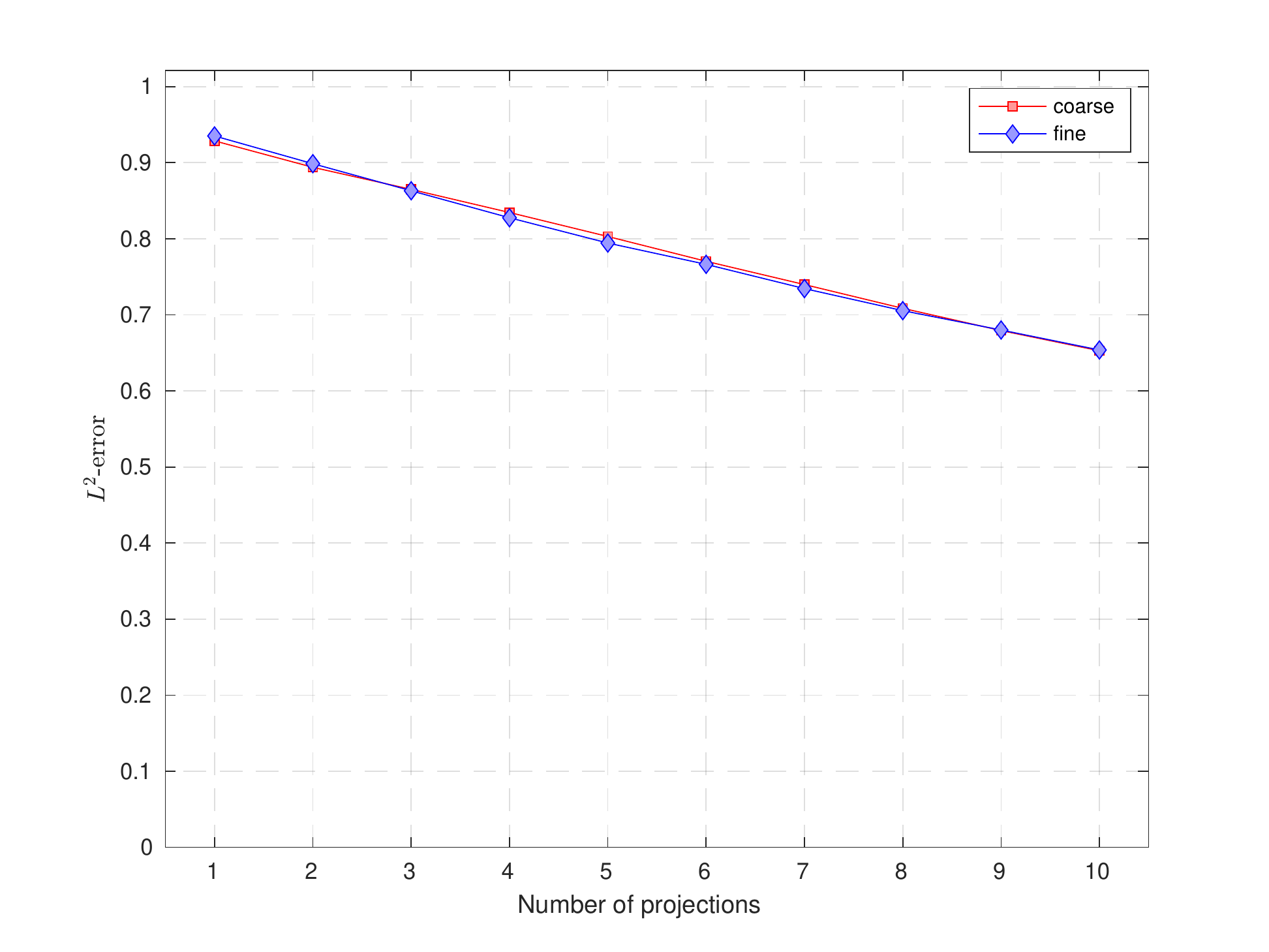}
  \end{center}
  \caption{{\sc Test~1.3.} Left: Obstruction.  Right: Mean $L^2(D)$ reconstruction errors over a sample of 1000 absorption targets as functions of the number of (noisy) projections. The diamonds correspond to the A-optimal projections deduced using the dense discretization and the squares to those deduced using the coarse discretization. The reconstruction and corresponding errors were computed using the dense discretization.}
      \label{fig:test1_3b}
\end{figure}

Next we test Algorithm~\ref{alg:basic_optimization} with an obstruction inside $D$ and also examine the effect that a coarse discretization has on the optimal angles. Only A-optimality is considered --- the results for D-optimality would be qualitatively the same in this case. The obstruction blocking the X-rays is defined via $D_{\rm obst} = \{ x \in D \ | \ x_1<0.5,\ 0.45 < x_2 < 0.55 \}$, with $D \setminus D_{\rm obst}$ being the ROI; see the left-hand image in Figure~\ref{fig:test1_3b}. The parallel beam source-receiver pair has width $0.5$, and the parameters defining the prior and the additive measurement noise are also the same as in the previous test: $\sigma = 0.02$, $\gamma = 1$ and $\ell = 0.05$. We employ two different levels of discretization: moderately dense, with $N=100$ and $m=23$, and very coarse, with $N=25$ and $m=6$.

\begin{figure}[ht!]
  \begin{center}
		\includegraphics[width = \columnwidth]{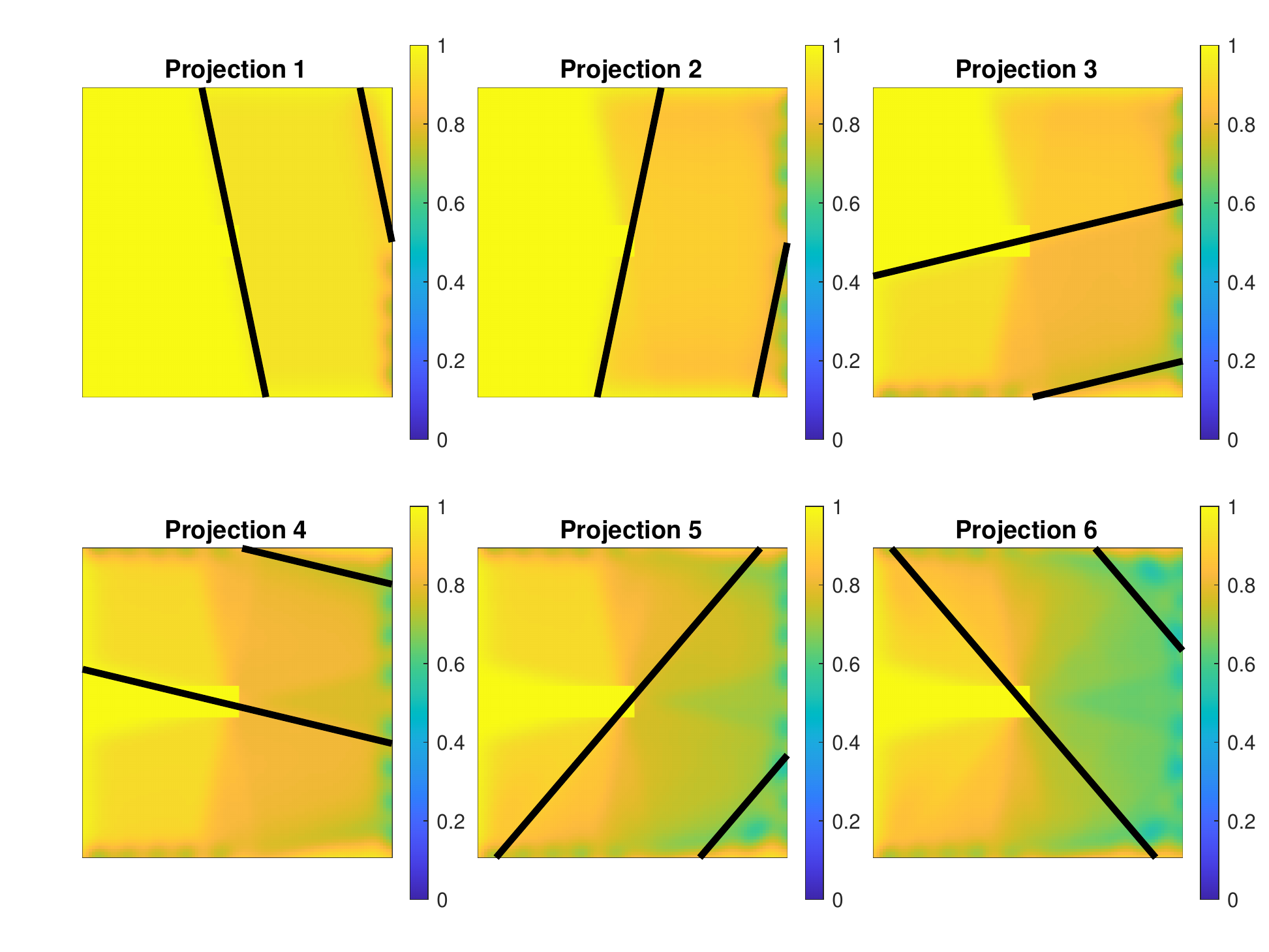}
  \end{center}
  \caption{{\sc Test~1.3.} The first six sequentially A-optimal parallel beam projections for the dense discretization with the corresponding pixelwise (posterior) standard deviations in the background.}
      \label{fig:test1_3a}
  \end{figure}

Figure~\ref{fig:test1_3a} presents the first six sequentially A-optimal projections for the denser discretization, with the corresponding pixelwise (posterior) standard deviations shown in the background. As expected, the algorithm tries to avoid imaging through the obstruction and emphasizes to begin with reducing the uncertainty in the right half of $D$, where imaging from almost all directions is possible. The projections produced by Algorithm~\ref{alg:basic_optimization} for the coarser discretization, not illustrated here, are qualitatively similar, but not exactly the same as the ones shown in Figure~\ref{fig:test1_3a}.

The right-hand image of Figure~\ref{fig:test1_3b} depicts the mean $L^2(D)$ reconstruction errors after each of the first ten (noisy) A-optimal projections over a sample of $1000$ target absorptions drawn from the assumed prior. Both sets of optimal projections, i.e.~those deduced using the dense discretization ($n=100$ and $m=23$) and the ones deduced using the coarse discretization ($n=25$ and $m=6$), are considered. However, all reconstructions and the corresponding $L^2(D)$ errors are computed using the dense discretization. Based on the right-hand image of Figure~\ref{fig:test1_3b}, it seems quite obvious that the level of discretization used for finding the sequentially A-optimal projections does not play a major role in the performance of the overall approach --- at least for the studied parameter values and the simple geometry.

\subsection{Test~2: Adaptive region of interest}
Let us then study a setting where the information on the location of the ROI is updated during the imaging process, that is, we test Algorithm~\ref{alg:ROI_optimization}. We adopt the geometry in the left-hand image of Figure~\ref{fig:test1_3b} as well as the parameter values in (the denser discretization of) Test~1.3: $\gamma=1$, $\ell = 0.05$, $\sigma = 0.02$, $N=100$ and $m=23$. The width of the parallel beam source-receiver pair is $0.5$, and $D\setminus D_{\rm obst}$ is the ROI to begin with.

\begin{figure}
  \begin{center}
		\includegraphics[width = \columnwidth]{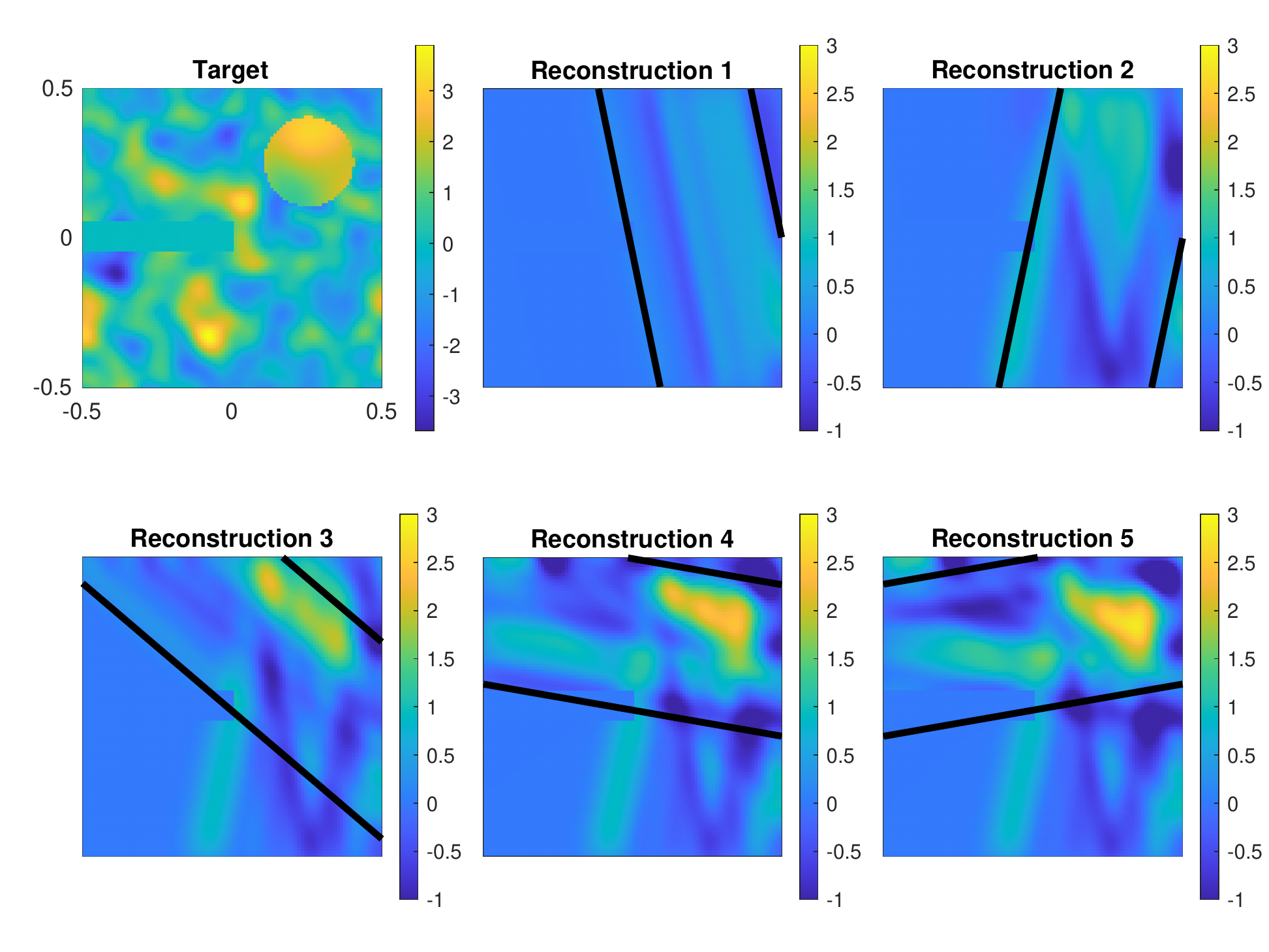}
  \end{center}
  \caption{{\sc Test~2.} Top left: True absorption and the obstruction $D_{\rm obst}$. Rest of the top row: The first two sequentially A-optimal projections and the resulting reconstructions from noisy measurements. Bottom row: The next three sequentially A-optimal projection after redefining the ROI to be the top right quadrant of $D$. The  resulting reconstructions from noisy measurements are shown in the background.}
      \label{fig:test2}
  \end{figure}

The top left image in Figure~\ref{fig:test2} shows the true absorption inside $D$:  in the top right corner, there is an inhomogeneity that is characterized by a higher level of absorption compared to the rest of the domain. The other two images in the top row of Figure~\ref{fig:test2} present the first two A-optimal projections and the corresponding reconstructions, i.e.~the posterior means, computed from simulated noisy measurements. Notice that these first projections are the same as the corresponding ones in Figure~\ref{fig:test1_3b}. Since the second reconstruction hints there may be something interesting going on in the top right quadrant of $D$, the ROI is redefined accordingly. The subsequent three A-optimal projections are illustrated in the bottom row of Figure~\ref{fig:test2} together with the resulting reconstructions. Comparing these results with those in Figure~\ref{fig:test1_3b}, it seems that the redefinition of the ROI had the desired effect as all of the last three projections clearly provide information on the top right quadrant of $D$.

\subsection{Test~3: Simultaneous estimation of the prior correlation length}
Our final numerical experiment applies Algorithm~\ref{alg:data-driven_adaptive_estimation} to deducing the value of an unknown correlation length $\rho = \ell$ in the prior covariance \eqref{eq:prior_cov}. The other free parameters are chosen to be $\gamma = 1$, $\sigma = 0.05$, $N=75$ and $m=39$. The parallel beam source-receiver pair has width $1$, i.e.~its position cannot be adjusted laterally, there are no obstructions inside $D$, and only A-optimality is considered. The presented numerical results would be qualitatively similar if the unknown parameter in the prior were instead the pointwise standard deviation or/and if D-optimality were used as the criterion for choosing the projection angles.

The following simple test is repeated $1000$ times: A random value for the true correlation length $\ell$ is picked from the uniform distribution over $[0.04,0.06]$ and a corresponding target absorption is subsequently drawn from the zero-mean Gaussian density with the prior covariance \eqref{eq:prior_cov}. Then, Algorithm~\ref{alg:data-driven_adaptive_estimation} is run with the conservative initial guess $\ell_0=0.15$ for the correlation length. The estimate for the correlation length after the first projection, i.e.~$\ell_1$, is determined by performing ten steps of the golden section line search over the interval $[0.01, 0.2]$ and subsequently applying the Newton's method to fine-tune the location of the optimum. When determining $\ell_k$ for $k \geq 2$, mere Newton's method is used with $\ell_{k-1}$ as the initial guess. In all cases, the Newton's method is stopped when its step size is less than $10^{-4}$, and as a safety measure, Newton steps with absolute value larger than $0.01$ are scaled down to $0.01$ (retaining the sign of the step). No obvious problems with the convergence of the implemented minimization routine for finding $\ell_k$ were observed during computations.

\begin{figure}
  \begin{center}
  \includegraphics[width = 0.45\columnwidth]{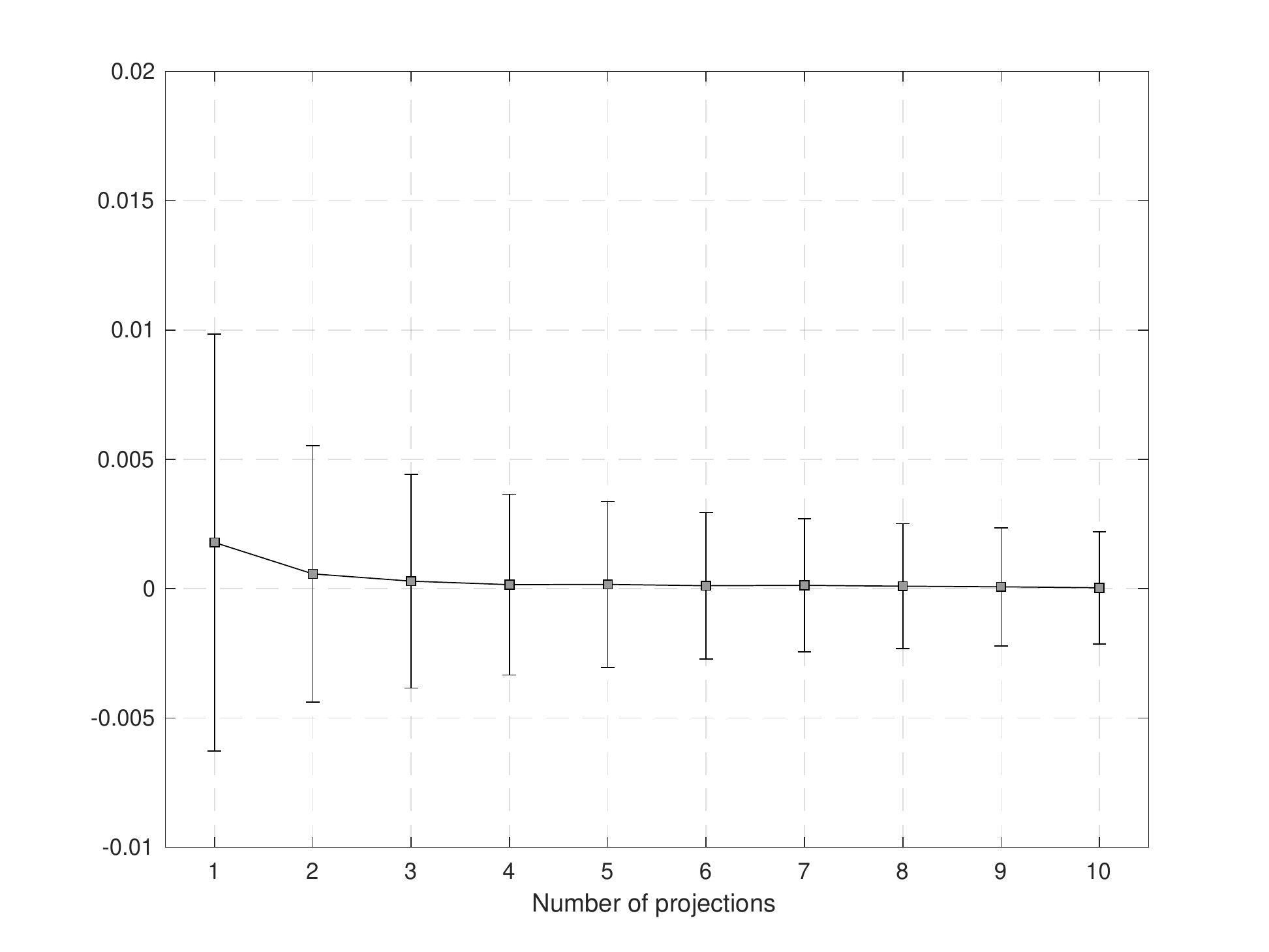}
  \includegraphics[width = 0.45\columnwidth]{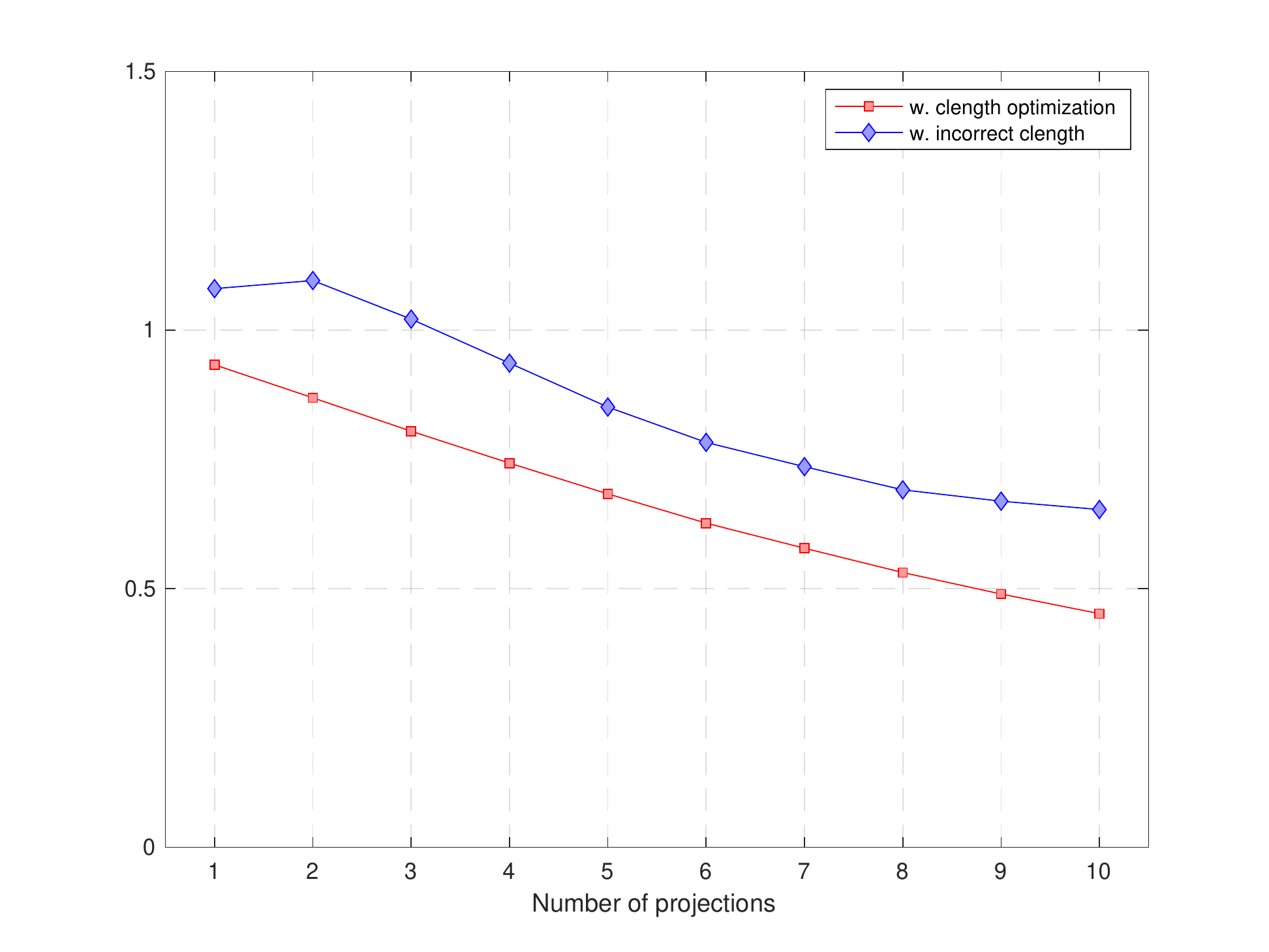}
  \end{center}
  \caption{{\sc Test~3.} Left: Average signed error in the estimate for the prior correlation length over a sample of $1000$ target absorptions as a function of the number of sequentially A-optimal projection angles. The errorbars indicate the intervals of two standard deviations computed over the studied sample. Right: Mean $L^2(D)$ reconstruction errors when the prior correlation length is simultaneously estimated by Algorithm~\ref{alg:data-driven_adaptive_estimation} (squares) and the corresponding mean $L^2(D)$ errors if the conservative initial guess $l_0 = 0.15$ is used both for optimizing the projection angles and for forming the reconstructions (diamonds).}
      \label{fig:test3}
\end{figure}

The left-hand imaged in Figure~\ref{fig:test3} shows the mean signed errors in the estimate for $\ell$ during the first ten iterations of Algorithm~\ref{alg:data-driven_adaptive_estimation} over the aforementioned sample of 1000 absorption targets. The corresponding standard deviations are visualized as errorbars. The mean signed error is almost zero independently of the number of (sequentially A-optimal) projections, whereas the standard deviation of the estimate over the sample of absorption targets decreases from approximately $0.008$ after the first projection to about $0.002$ after the tenth one. Altogether, the deduction of the prior correlation length based on noisy projection data seems to function well.

Certain mean $L^2(D)$ errors in the reconstructions computed based on noisy parallel beam projections of the studied sample of target absorptions are shown as functions of the number of projections in the right-hand image of Figure~\ref{fig:test3}. The reconstructions and the corresponding $L^2(D)$ errors are computed in two different ways: (i) by accounting for the change in the estimate for the prior correlation length predicted by Algorithm~\ref{alg:data-driven_adaptive_estimation} after each new sequentially A-optimal projection and (ii) by sticking with the conservative initial guess of $0.15$ for the correlation length throughout the optimization of the projection angles and in the computation of the corresponding reconstructions. It is obvious that including the estimation of the correlation length as a part of the reconstruction algorithm leads to superior reconstruction accuracy, although this also makes the algorithm significantly more expensive computationally. If the initial guess $\ell_0$ for the prior correlation length were (significantly) more accurate, it is debatable whether fine-tuning the estimate for $\ell$ would be worthwhile computationally.

\section{Concluding remarks}
\label{sec:concl}
We have introduced a greedy, exhaustive optimization algorithm that applies Bayesian OED to optimizing parallel beam projections in X-ray imaging. Although the sequential algorithm does not guarantee global optimality for the set of projections it predicts, our numerical experiments demonstrate that it usually exhibits satisfactory global behavior as well. More importantly, it fits to the medical imaging paradigm where the patient is exposed to the minimal amount of radiation without knowing {\em a priori} the necessary number of imaging angles. To this end, we have discussed modifications of the algorithm that allow tuning the ROI or some parameters in the (original) prior based on the measured data.

The main obstacle for adopting our algorithm in practical use is arguably its computational complexity: as the optimization target typically has many local minima (cf.~Figure~\ref{fig:test1_1a}), we seek the optimal parameters via evaluating the target on a dense enough grid in order to guarantee finding the global optimum for the next projection. In a three-dimensional setting, such a simple approach becomes prohibitively expensive as the number of pixels in a single projection image can be in the order of $m=10^5$, and thus a more sophisticated adaptive optimization routine needs to be developed for more realistic imaging geometries. On the positive side, the overall performance of our approach does not seem to depend heavily on the level of discretization, which may permit a procedure where the next projection geometry is optimized using a coarse discretization, but the actual reconstruction corresponds to a much denser one.

In addition, only a very simple Gaussian prior for the absorption distribution was employed in our numerical studies. The introduced algorithm should thus be tested with more realistic priors and, even more importantly, its functionality should be verified if no Gaussian prior is well aligned with the imaged target.

\appendix

\section{Target functionals for Bayesian experimental design}
\label{sec:KL}
The purpose of this appendix is to deduce the minimization targets \eqref{eq:Aoptimal} and \eqref{eq:Doptimal} and to indicate how the latter is actually related to the information gain when a prior is replaced by the posterior. It is also explained why $\Gamma_{\rm post}$ should be replaced by $(\Gamma_{\rm post})_{11}= \Gamma_{\rm ROI}$ in \eqref{eq:Doptimal} if one is only interested in the information gain over the ROI,~cf.~\eqref{eq:block_form}.

\subsection{D-optimality}
The Kullback--Leibler divergence of the distribution of $X$ from that of $X \,| \, y$,~i.e.,
\begin{align}
  \label{eq:KL}
  D_{\rm KL}(X | \, y \, \| \, X; p) &= \int_{\R^n} \log \left(\frac{\pi(x\, |\, y;p)}{\pi_{\rm pr} (x)}\right)  \pi(x\, |\, y;p) \,\dx \nonumber \\
  &= \int_{\R^n} \big( \log (\pi(x\, |\, y;p)) - \log (\pi_{\rm pr}(x)) \big) \, \pi(x\, |\, y;p)  \, \dx ,
\end{align}
has an interpretation as the increase in the level of information when the prior distribution of $X$ is replaced by the posterior induced by the data $y$. The aim of D-optimal design is to maximize the expectation of $D_{\rm KL}(X | \, y \, \| \, X; p)$ over all measurements, which means the precise form of the {\em minimization} target is
$$
\tilde{\Phi}_D(p) = \mathbb{E}^y\big( - D_{\rm KL}(X | \, y \, \| \, X; p) \big) = - \int_{\R^m} D_{\rm KL}(X | \, y \, \| \, X; p)\, \pi(y; p) \, \dy;
$$
see~\eqref{eq:Dtarget} for comparison.

The expression
\begin{equation}
  \label{eq:diff_entropy}
 h( X \, | \, y) = -\int_{\R^n} \log (\pi(x \, | \, y\, ; p)) \, \pi (x \, | \, y\, ;p)\, \dx
\end{equation}
appearing in \eqref{eq:KL} is called the \textit{differential entropy} of $X \, | \, y$. Since $X \, | \, y$ is Gaussian, $h( X \, | \, y)$ allows an explicit representation~\cite{Shannon48}
\begin{equation}
  \label{eq:post_simple}
  h( X \, | \, y)
  = \frac{n}{2} + \frac{n}{2}\log(2\pi) + \frac{1}{2}\log\big(\!\det(\Gamma_{\rm post}(p)) \big).
\end{equation}
On the other hand, the expected value of the second term on the right-hand side of \eqref{eq:KL} satisfies
\begin{align*}
\mathbb{E}^y \Big(\int_{\R^n} \log (\pi_{\rm pr} (x)) \, \pi(x\, |\, y ; p)\, \dx\Big) = &\int_{\R^m} \int_{\R^n} \log (\pi_{\rm pr} (x)) \, \pi(x, y ;p)\, \dx\, \dy\\[1mm]
  = &\int_{\R^n} \int_{\R^m} \pi(y \, | \, x; p) \, \dy \, \log (\pi_{\rm pr} (x) )\,\pi_{\rm pr} (x) \, \dx\\[1mm]
  = & \int_{\R^n} \log (\pi_{\rm pr} (x)) \pi_{\rm pr}(x)\, \dx\\[1mm]
  = &-\frac{n}{2}-\frac{n}{2}\log (2\pi) -\frac{1}{2}\log(\det \Gamma_{\rm pr}),
\end{align*}
where the last expression is simply $-h(X)$ for the Gaussian random variable $X$. In particular, take note that this expression does not depend on the design variable $p$.

Because $h( X \, | \, y)$ does not depend on $y$, we finally arrive at the simple expression
\begin{align}
  \label{eq:Dtarget2}
\tilde{\Phi}_D(p)  & =  \mathbb{E}^y \big(h( X \, | \, y)\big)  - h(X)  \\[1mm]
  \nonumber &=  \frac{1}{2}\big(\log\big(\! \det \Gamma_{\rm post}(p)\big) - \log(\det\Gamma_{\rm pr}) \big) ,
\end{align}
which explains why the minimizers of $\tilde{\Phi}_D(p)$ and $\Phi_D(y)$ in \eqref{eq:Dtarget} are the same if ROI is the whole  $D$. Note that $-\tilde{\Phi}_D(p)$ corresponds to the actual information gain that is to be maximized.

Finally, observe that the above calculation can be repeated as such if $\pi_{\rm pr}(x)$ and $\pi(x \, | \, y; p)$ are replaced by the corresponding marginal densities over the ROI. Hence, the precise form of the minimization target when the information gain is to be maximized over the ROI is obtained by replacing $\Gamma_{\rm post}(p)$ by $(\Gamma_{\rm post}(p))_{11}$ and $\Gamma_{\rm pr}$ by $(\Gamma_{\rm pr})_{11}$ in \eqref{eq:Dtarget2}, cf.~\eqref{eq:Dtarget}.

\subsection{A-optimality}
To begin with, let $A \in \R^{l \times n}$ and observe that the CM estimate $\widehat{x}(p) = \widehat{x}(y; p) \in \R^n$, given in \eqref{eq:posterior_cov}, naturally depends on the data $y \in \R^m$. The expected squared distance of the unknown $X$ from $\widehat{x}(Y; p)$ in the squared seminorm induced by $A^T \! A$ over the joint distribution of $X$ and $Y$ reads as follows:
\begin{align*}
  \Phi_{\rm A}(p) &= \int_{\R^m} \int_{\R^n} \big(x - \widehat{x}(y; p)\big)^T \! A^T \!A \big(x - \widehat{x}(y; p) \big) \, \pi(x,y; p) \, \dx \, \dy \\[1mm]
  &= \int_{\R^m} \int_{\R^n} {\rm tr}\big((x - \widehat{x}(y; p))^T \! A^T \!A (x - \widehat{x}(y; p) )\big) \, \pi(x \, | \, y; p) \, \dx \, \pi(y; p) \, \dy
\end{align*}
as the trace of a scalar is the scalar itself. Because the trace is linear and invariant under cyclic permutations, we get
\begin{align*}
  \Phi_{\rm A}(p) &= {\rm tr}\Big(A \int_{\R^m} \int_{\R^n} (x - \widehat{x}(y; p)) (x - \widehat{x}(y; p))^T \, \pi(x \, | \, y; p) \, \dx \, \pi(y; p) \, \dy \,  A^T \Big) \\[1mm]
  &= {\rm tr}\Big(A \int_{\R^m} \Gamma_{\rm post}(p) \, \pi(y; p) \, \dy \,  A^T \Big) \\[1mm]
  &= {\rm tr} ( A \Gamma_{\rm post}(p) A^T),
\end{align*}
where the last step follows from $\Gamma_{\rm post}(p)$ being independent of $y$ in our simple, Gaussian and linear setting.

\bibliographystyle{acm}
\bibliography{optdir-refs}
\end{document}

%% file: figures_final/circle.tex
\centering
\begin{tikzpicture}[scale=2.0]
  \filldraw[fill=blue!20!white] (-1, -1) rectangle (1, 1);
  \draw[fill=green!20!white] (0.2, 0.2) circle (0.5) node {ROI};
\end{tikzpicture}

%% file: figures_final/one_obstruction.tex
\centering
\begin{tikzpicture}[scale=2]
  \filldraw[fill=green!20!white] (-1, -1) rectangle (1, 1);
  \filldraw[fill=red!20!white] (-1, -0.1) rectangle (0, 0.1);
  \node[] at (0, 0.5) {ROI};
  \node[] at (-0.5, 0) {Obs.};
\end{tikzpicture}